\documentclass[11pt]{article}
\usepackage{latexsym,amssymb,mathrsfs}

\usepackage{psfrag}
\usepackage{amsmath}
\usepackage{color}

\newtheorem{Theorem}{\bf Theorem}[section]
\newtheorem{Lemma}{\bf Lemma}[section]
\newtheorem{Proposition}{\bf Proposition}[section]
\newtheorem{Corollary}{\bf Corollary}[section]
\newtheorem{Remark}{\bf Remark}[section]
\newtheorem{Example}{\bf Example}[section]
\newtheorem{Definition}{\bf Definition}[section]

\newenvironment{theorem}{\begin{Theorem}$\!\!\!$}{\end{Theorem}}
\newenvironment{lemma}{\begin{Lemma}$\!\!\!$}{\end{Lemma}}
\newenvironment{proposition}{\begin{Proposition}$\!\!\!$}{\end{Proposition}}

\newenvironment{remark}{\begin{Remark}$\!\!\!$}{\end{Remark}}

\newenvironment{definition}{\begin{Definition}$\!\!\!$}{\end{Definition}}

\numberwithin{equation}{section}

\def\XXint#1#2#3{{\setbox0=\hbox{$#1{#2#3}{\int}$}
\vcenter{\hbox{$#2#3$}}\kern-.5\wd0}}

\begin{document}

\title{The Cauchy problem for the Finsler heat equation}

\author{Goro Akagi, Kazuhiro Ishige and Ryuichi Sato}
\date{}
\maketitle
\begin{abstract}
Let $H$ be a norm of ${\bf R}^N$ and $H_0$ the dual norm of $H$. 
Denote by $\Delta_H$ the Finsler-Laplace operator defined by $\Delta_Hu:=\mbox{div}\,(H(\nabla u)\nabla_\xi H(\nabla u))$. 
In this paper we prove that 
the Finsler-Laplace operator $\Delta_H$ acts as a linear operator to $H_0$-radially symmetric smooth functions. 
Furthermore, we obtain an optimal sufficient condition for the existence of the solution 
to the Cauchy problem  for  the Finsler heat equation
$$
\partial_t u=\Delta_H u,\qquad x\in{\bf R}^N,\quad t>0,
$$
where $N\ge 1$ and $\partial_t:=\partial/\partial t$. 
\end{abstract}
\vspace{15pt}
\noindent Addresses:

\smallskip
\noindent G. A.:  Mathematical Institute, Tohoku University,
Aoba, Sendai 980-8578, Japan\\
\noindent 
E-mail: {\tt akagi@m.tohoku.ac.jp}\\

\noindent K. I.:  Mathematical Institute, Tohoku University,
Aoba, Sendai 980-8578, Japan\\
\noindent 
E-mail: {\tt ishige@m.tohoku.ac.jp}\\

\noindent R. S.: Mathematical Institute, Tohoku University,
Aoba, Sendai 980-8578, Japan\\
\noindent 
E-mail: {\tt ryuichi.sato.d4@tohoku.ac.jp}\\
\vspace{15pt}
\vspace{3pt}
\newpage
\section{Introduction}
Let $N\ge 1$ and let $H\in C({\bf R}^N)\cap C^1({\bf R}^N\setminus\{0\})$ be a norm of ${\bf R}^N$, that is 
\begin{equation}
\label{eq:1.1}
\left\{
\begin{array}{l}
 \mbox{$H\ge 0$ in ${\bf R}^N$ and $H(\xi)=0$ if and only if $\xi=0$},\vspace{3pt}\\
 \mbox{$H$ is convex in ${\bf R}^N$},\vspace{3pt}\\
 \mbox{$H(\alpha\xi)=|\alpha|H(\xi)$ for $\xi\in{\bf R}^N$ and $\alpha\in{\bf R}$}.  
\end{array}
\right.
\end{equation}
We denote by $H_0$ the dual norm  of $H$ defined by
$$
H_0(x):=\sup_{\xi\in{\bf R}^N\setminus\{0\}}\frac{x\cdot\xi}{H(\xi)}. 
$$
Then 
\begin{equation}
\label{eq:1.2}
|x\cdot\xi|\le H_0(x)H(\xi),
\qquad
H(\xi)=\sup_{x\in{\bf R}^N\setminus\{0\}}\frac{x\cdot\xi}{H_0(x)}. 
\end{equation}
For any $x\in{\bf R}^N$, $\xi\in{\bf R}^N$ and $R>0$, we set 
\begin{equation*}
\begin{split}
B_H(\xi,R):= & \{\eta\in{\bf R}^N\,:\,H(\eta-\xi)<R\},\\
B_{H_0}(x,R):= & \{y\in{\bf R}^N\,:\,H_0(y-x)<R\}.
\end{split}
\end{equation*}
Throughout this paper we assume that 
\begin{equation}
\label{eq:1.3}
\mbox{$B_H(0,1)$ is strictly convex},
\end{equation}
which is equivalent to $H_0\in C^1({\bf R}^N\setminus\{0\})$ 
(see \cite[Corollary~1.7.3]{S}). 
\vspace{3pt}

Let $\Delta_H$ be the Finsler-Laplace operator associated with the norm~$H$, 
that is  
$$
\Delta_H\,u:=\mbox{div}(\nabla_\xi V(\nabla u))=\mbox{div}(H(\nabla u)\nabla_\xi H(\nabla u)),
$$
where $V(\xi):=H(\xi)^2/2$. 
We remark that $\Delta_H=\Delta$ if $H(\xi)=|\xi|$ for $\xi\in{\bf R}^N$. 
The Finsler-Laplace operator has been treated by many mathematicians
from various points of view 
(see e.g., \cite{BCS}, \cite{BP}, \cite{BGS}, \cite{CS}, \cite{CFV}, \cite{DPB}, \cite{FK}, \cite{OS1}, \cite{OS2}, \cite{T}, \cite{X} 
and references therein). 
This paper is concerned with the Cauchy problem for 
the Finsler heat equation
\begin{equation}
\label{eq:1.4}
\partial_t u=\Delta_H\,u,\qquad x\in{\bf R}^N,\quad t>0, 
\end{equation}
which is introduced as a gradient flow of the energy 
\begin{equation*}
I[u]:=\frac{1}{2}\int_{{\bf R}^N}H(\nabla u)^2\,dx.  
\end{equation*}
Equation~\eqref{eq:1.4} is a nonlinear parabolic equation with the following nice property: 
\begin{equation}
\label{eq:1.5}
\begin{split}
 & \mbox{If $u$ is a solution of \eqref{eq:1.4}, then $ku$ and $u(kx,k^2 t)$ are also solutions of \eqref{eq:1.4}}\\
 & \mbox{for any $k\in{\bf R}$}
\end{split}
\end{equation}
(see Definition~\ref{Definition:1.1} and Section~2).

\vspace{5pt}

 The Finsler-Laplace operator enjoys nice properties for ``smooth radial'' functions. To be more precise, 
for any function $v$ in ${\bf R}^N$, 
we say that $v$ is \emph{$H_0$-radially symmetric}  in ${\bf R}^N$ if 
there exists a function $v^\sharp$ on $[0,\infty)$ such that  
$$
v(x)=v^\sharp(r)\quad\mbox{for $x\in{\bf R}^N$ with $r=H_0(x)$}.
$$ 
Set $v^*(x):=v^\sharp(|x|)$ for $x\in{\bf R}^N$. 
For $k\in\{0,1,2,\dots\}$, 
we say that an $H_0$-radially symmetric function $v$ is  \emph{$C^k_{H_0}$-smooth}  
if $v^*\in C^k({\bf R}^N)$.  
One of  main  purposes of this paper is to obtain the following two nice properties of $H_0$-radially symmetric 
$C^2_{H_0}$-smooth functions. 
\begin{theorem}
\label{Theorem:1.1}
Assume \eqref{eq:1.1} and \eqref{eq:1.3}. 
Let $v$ and $w$ be $H_0$-radially symmetric $C^2_{H_0}$-smooth functions in ${\bf R}^N$. 
Then 
\begin{itemize}
  \item[{\rm (a)}] 
  $(\Delta_H v)(x)=(\Delta v^*)(y)$ for $x\in{\bf R}^N$ with $y=H_0(x)$;
  \item[{\rm (b)}] 
  $\Delta_H(\alpha v+\beta w)=\alpha\Delta_H v+\beta\Delta_H w$ in ${\bf R}^N$ 
  for any $\alpha$, $\beta\in{\bf R}$. 
\end{itemize}
\end{theorem}
The Finsler-Laplace operator is a nonlinear operator and often has a  complicated  form. 
Nevertheless, Theorem~\ref{Theorem:1.1} implies that 
the Finsler-Laplace operator acts as a linear operator to $H_0$-radially symmetric $C^2_{H_0}$-smooth functions. 
Furthermore, Theorem~\ref{Theorem:1.1} enables us to 
bring rich mathematical results for radially symmetric functions  to the study of  
$H_0$-radially symmetric functions. 
See also Section~2. 
\vspace{3pt}

The other purpose of this paper is to obtain an optimal growth condition on initial data 
for the existence of  solutions to  the Cauchy problem for the Finsler heat equation
\begin{equation}
\label{eq:1.6}
\partial_t u=\Delta_H\,u, \quad x\in{\bf R}^N,\,\,t>0,
\qquad
u(\cdot,0)=\mu\quad\mbox{in}\quad{\bf R}^N,
\end{equation}
where $\mu$ is a (signed) Radon measure in ${\bf R}^N$. 
The study of the growth conditions on initial data 
for the existence of solutions to parabolic equations is a classical subject. 
The growth conditions generally depend on the diffusion and the nonlinear terms in the parabolic problems 
(see e.g., \cite{AD1}--\cite{Ar}, \cite{BaP}, \cite{BCP}, \cite{DH1}, \cite{DH2}, \cite{GU}--\cite{IJK}, \cite{IS2}, \cite{W} and \cite{Y}). 
For the heat equation, the following holds. 
\begin{itemize}
  \item[{\rm (H1)}]
  Let $u$ be a nonnegative solution of $\partial_t u=\Delta u$ in ${\bf R}^N\times(0,T)$, where $T>0$. 
  Then there exists a unique nonnegative Radon measure $\mu$ in ${\bf R}^N$ such that 
  $$
  \lim_{t\to +0}\int_{{\bf R}^N}u(y,t)\phi(y)\,dy=\int_{{\bf R}^N}\phi(y)\,d\mu(y),
  \qquad\phi\in C_0({\bf R}^N). 
  $$
  Furthermore, $\mu$ satisfies 
  $$
  \sup_{x\in{\bf R}^N}\int_{B(x,1/\sqrt{\Lambda})}e^{-\Lambda|x|^2}d\mu<\infty\quad\mbox{for some $\Lambda>0$}.
  $$
  \item[{\rm (H2)}] 
  Let $\mu$ be a (signed) Radon measure in ${\bf R}^N$ such that 
  \begin{equation}
  \label{eq:1.7}
  \sup_{x\in{\bf R}^N}\int_{B(x,1/\sqrt{\Lambda})}e^{-\Lambda|x|^2}d|\mu|<\infty\quad\mbox{for some $\Lambda>0$}. 
  \end{equation}
  Then 
  $$
  u(x,t):=(4\pi t)^{-\frac{N}{2}}\int_{{\bf R}^N}\exp\left(-\frac{|x-y|^2}{4t}\right)\,d\mu(y)
  $$
  satisfies
  \begin{equation}
  \label{eq:1.8}
  \partial_t u=\Delta u\quad\mbox{in}\quad{\bf R}^N\times(0,1/4\Lambda),
  \qquad
  u(\cdot,0)=\mu\quad\mbox{in}\quad{\bf R}^N. 
  \end{equation}
  \item[{\rm (H3)}] 
  Let $\Lambda>0$ and set
  $$
  v(x,t):=(1-4\Lambda t)^{-\frac{N}{2}}
  \exp\left(\frac{\Lambda|x|^2}{1-4\Lambda t}\right). 
  $$
  Then $v$ is a solution of \eqref{eq:1.8} and it satisfies 
  \eqref{eq:1.7} and 
  $$
  \min_{x\in{\bf R}^N}v(x,t)=v(0,t)\to\infty\quad\mbox{as}\quad t\to\frac{1}{4\Lambda}.
  $$ 
\end{itemize}
See \cite{Ar} and \cite{W}. (See also \cite{I} and  \cite{IJK}).

Consider the Cauchy problem for the porous medium equation 
\begin{equation}
\label{eq:1.9}
\partial_t u=\Delta(|u|^{m-1}u),\quad x\in{\bf R}^N,\,\,t>0,
\qquad
u(\cdot,0)=\mu\quad\mbox{in}\quad{\bf R}^N, 
\end{equation}
where $m>1$. 
B\'enilan, Crandall and Pierre \cite{BCP} proved that 
problem~\eqref{eq:1.9} possesses a solution provided that 
$$
\limsup_{\rho\to\infty}\,\rho^{-N-\frac{2}{m-1}}\int_{B(0,\rho)}d\,|\mu|<\infty. 
$$
On the other hand, in the case of $(N-2)_+/2<m<1$, 
Herrero and Pierre~\cite{HP} 
proved that problem~\eqref{eq:1.9} possesses a solution for any $L^1_{{\rm loc}}({\bf R}^N)$ initial data. 
For similar results to  parabolic  $p$-Laplace  equations  
and more general nonlinear parabolic equations, see e.g., 
\cite{AD1}, \cite{DH1}, \cite{DH2}, \cite{I0}, \cite{I} and \cite{I1}. 
\vspace{7pt}

We formulate  a  definition of  solutions to  \eqref{eq:1.6}. 
\begin{definition}
\label{Definition:1.1}
Let $u$ be a measurable function $u$ in ${\bf R}^N\times(0,T)$, where $T>0$. 
\begin{itemize}
  \item[{\rm (i)} ]
  We say that $u$ is a solution of \eqref{eq:1.4} in ${\bf R}^N\times(0,T)$ if 
  \begin{equation}
  \label{eq:1.10}
  u\in C((0,T):L^1(B_{H_0}(0,R)))\,\cap\,L^1((0,T-\epsilon):W^{1,1}(B_{H_0}(0,R)))
  \end{equation}
  for any $R>0$ and $0<\epsilon<T$
  and $u$ satisfies 
  \begin{equation*}
  \begin{split} 
   & \int_{{\bf R}^N}u(y,t)\varphi(y,t)\,dy
  +\int^t_\tau\int_{{\bf R}^N}\,
  \left[-u\partial_t\varphi+H(\nabla u)\nabla_\xi H(\nabla u)\nabla\varphi\right]\,dy\,ds\\
   & =\int_{{\bf R}^N}u(y,\tau)\varphi(y,\tau)\,dy 
  \end{split}
  \end{equation*}
  for all $0<\tau<t<T$ and $\varphi\in C^\infty_0({\bf R}^N\times(0,T))$.
  \item[{\rm (ii)} ]
  Let $\mu$ be a {\rm(}signed\,{\rm)} Radon measure in ${\bf R}^N$. 
  Then we say that $u$ is a solution of \eqref{eq:1.6} in ${\bf R}^N\times[0,T)$
  if $u$ satisfies \eqref{eq:1.10} and 
  \begin{equation*}
  \begin{split}
   & \int_{{\bf R}^N}u(y,t)\varphi(y,t)\,dy
  +\int_0^t\int_{{\bf R}^N}\,\left[-u\partial_t\varphi+H(\nabla u)\nabla_\xi H(\nabla u)\nabla\varphi\right]\,dy\,ds\\
   & =\int_{{\bf R}^N}\varphi(y,0)\,d\mu(y)
  \end{split}
  \end{equation*}
  for all $0<t<T$ and $\varphi\in C^\infty_0({\bf R}^N\times[0,T))$.
\end{itemize}
\end{definition}

Now we are ready to state the second theorem of this paper. 
In  the following  Theorem~\ref{Theorem:1.2} 
we obtain similar results as in statements~(H1), (H2) and (H3) for the Cauchy problem 
to the Finsler heat equation. 
\begin{theorem}
\label{Theorem:1.2}
Assume \eqref{eq:1.1} and \eqref{eq:1.3}. 
\begin{itemize}
  \item[{\rm (i)}] 
  Let $u$ be a nonnegative solution of \eqref{eq:1.4} in ${\bf R}^N\times(0,T)$ for some $T>0$. 
  Then there exists a unique nonnegative Radon measure $\mu$ in ${\bf R}^N$ such that
  \begin{equation}
  \label{eq:1.11}
  \lim_{t\to +0}\int_{{\bf R}^N}u(y,t)\phi(y)\,dy=\int_{{\bf R}^N}\phi(y)\,d\mu(y),
  \qquad\phi\in C_0({\bf R}^N). 
  \end{equation}
  Furthermore, 
  $$
  \sup_{x\in{\bf R}^N}\int_{B_{H_0}(x,1/\sqrt{\Lambda})} e^{-\Lambda H_0(y)^2}d\mu(y)<\infty
  \quad\mbox{for some $\Lambda>0$}.
  $$
  \item[{\rm (ii)}]
  Let $\mu$ be a {\rm(}signed\,{\rm)} Radon measure in ${\bf R}^N$ such that 
  \begin{equation}
  \label{eq:1.12}
  \sup_{x\in {\bf R}^N}\int_{B_{H_0}(x,1/\sqrt{\Lambda})} e^{-\Lambda H_0(y)^2}\,d|\mu|(y)<\infty
  \quad\mbox{for some $\Lambda>0$}.
  \end{equation}
  Then there exists a nonnegative solution $u$ of \eqref{eq:1.6} in ${\bf R}^N\times[0,S_\Lambda)$, 
  where $S_\Lambda:=1/4\Lambda$, such that
  \begin{itemize}
  \item[{\rm (a)}]
  $u\in C^{1,\alpha;0,\alpha/2}({\bf R}^N\times(0,S_\Lambda))$ for some $\alpha \in (0,1)$\/{\rm ;}
  \item[{\rm (b)}] 
  For any $\lambda>\Lambda$, there exists a constant $C$ such that 
  \begin{equation}
  \label{eq:1.13}
   \sup_{0 < t < S_\lambda}\int_{{\bf R}^N} e^{-g_\lambda(x,t)}|u(y,t)|\,dy
  \le C\sup_{x\in {\bf R}^N}\int_{B_{H_0}(x,1/\sqrt{\Lambda})} e^{-\Lambda H_0(y)^2}\,d|\mu|(y).
  \end{equation}
	       Here $g_\lambda(x,t):=\lambda H_0(x)^2/(1-4\lambda t)$. 
  \end{itemize}
\end{itemize}
\end{theorem}
Theorem~\ref{Theorem:1.3} assures that $S_\Lambda$ is the optimal maximal existence time 
of the solution of \eqref{eq:1.6} under assumption~\eqref{eq:1.12}. 
\begin{theorem}
\label{Theorem:1.3}
Let $\Lambda>0$ and set 
$$
v(x,t):=(1-4\Lambda t)^{-\frac{N}{2}}
\exp\left(\frac{\Lambda H_0(x)^2}{1-4\Lambda t}\right).  
$$
Then $v$ is a solution of \eqref{eq:1.4} in ${\bf R}^N\times[0,S_\Lambda)$ 
and it satisfies \eqref{eq:1.12} and  
$$
\min_{x\in{\bf R}^N}v(x,t)=v(0,t)\to\infty\quad\mbox{as}\quad t\to S_\Lambda.
 $$
\end{theorem}

We explain  an  idea of the proof of Theorem~\ref{Theorem:1.2}.  
Since $\Delta_H$ is a nonlinear operator, 
we  cannot  apply the standard theory for linear parabolic equations to the Finsler heat equation. 
On the other hand, 
setting 
$$
A(\xi):=H(\xi)\nabla_\xi H(\xi)\quad\mbox{for}\quad \xi\not=0,
\qquad
A(\xi):=0\quad\mbox{for}\quad\xi=0,
$$
we have
\begin{equation}
\label{eq:1.14}
\begin{split}
 &  A \in C({\bf R}^N ; {\bf R}^N),  \quad \Delta_H u=\mbox{div}\,A(\nabla u) \ \mbox{ in } \mathscr{D}', \\
 & A(\xi)\cdot \xi=H(\xi)^2\quad\mbox{and}\quad H_0(A(\xi))=H(\xi)\quad\mbox{for}\quad\xi\in{\bf R}^N. 
 \end{split}
 \end{equation}
Since $H$ and $H_0$ are equivalent to the  Euclidean  norm of ${\bf R}^N$, 
it follows from \eqref{eq:1.14} that
\begin{equation}
\label{eq:1.15}
A(\xi)\cdot\xi\ge C_1|\xi|^2\quad\mbox{and}\quad |A(\xi)|\le C_2|\xi|\quad\mbox{for}\quad\xi\in{\bf R}^N,
\end{equation}
where $C_1$ and $C_2$ are positive constants. 
Then we can apply the arguments in \cite[Theorem~1.8]{IJK} 
to obtain assertion~(i). 

For the proof of assertion~(ii), 
by property~\eqref{eq:1.5} it suffices to consider the case of $\Lambda<1$. 
We introduce some new techniques to obtain uniform estimates of approximate solutions. 
Firstly, we construct approximate solutions $\{u_n\}$ by using  a subdifferential formulation, 
in order to preserve structure of the Finsler-Laplace operator (see, e.g., \eqref{eq:1.5}). 
Secondly, for any sufficiently small $\delta>0$, 
we obtain uniformly local $L^{1+\delta}$ estimates of $\{u_n\}$ 
with the aid of the Besicovitch covering theorem 
(see Lemma~\ref{Lemma:3.2}). 
Then we  modify  the arguments  used  in \cite{I0}--\cite{I1} to obtain uniformly local $L^1$ estimates of 
the functions 
$$
U_n:=F\left(e^{-H_0(y)^2(1+s^\ell)}|u_n(y,s)|\right),
$$
where $0 < \ell < 1$ and 
$$
F(v):=\int_0^v \min\{\tau^\delta,1\}\,d\tau
=\left\{
\begin{array}{ll}
v^{1+\delta}/(1+\delta) & \mbox{for}\quad 0\le v\le 1,\vspace{3pt}\\
v-\delta/(1+\delta) & \mbox{for}\quad v>1.
\end{array}
\right.
$$
Furthermore,  employing  property~\eqref{eq:1.5},  
we improve uniformly local $L^1$ estimates  for  $U_n$, which also asserts those for  
the functions $e^{-H_0(y)^2(1 + s^\ell)}u_n$. 
Finally, carrying out  a similar argument as in \cite{I0}--\cite{I1}, 
we show the existence of solutions  to  \eqref{eq:1.6} and complete the proof of assertion~(ii).

The rest of this paper is organized as follows. 
In Section~2 we recall some basic properties of $H$ and $H_0$ 
and prove Theorems~\ref{Theorem:1.1} and \ref{Theorem:1.3}. 
In Section~3 we obtain a priori estimates of the solutions of \eqref{eq:1.6} 
by improving the arguments in \cite{I}. 
In Section~4 we complete the proof of Theorem~\ref{Theorem:1.2}. 
\section{Proof of Theorems~\ref{Theorem:1.1} and \ref{Theorem:1.3}}
In this section we recall some basic properties of $H$ and $H_0$ and prove Theorems~\ref{Theorem:1.1} and \ref{Theorem:1.3}. 
Assume \eqref{eq:1.1} and \eqref{eq:1.3}. Then $H$, $H_0\in C^1({\bf R}^N\setminus\{0\})$.  
Furthermore, we have
\begin{equation}
\label{eq:2.1}
\left\{
\begin{array}{ll}
\xi\cdot\nabla_\xi H(\xi)=H(\xi), & \xi\in{\bf R}^N,\vspace{3pt}\\
\nabla_\xi H(t\xi)=\mbox{sign}\,(t)\,\nabla_\xi H(\xi), & \xi\in{\bf R}^N\setminus\{0\},\,\,\,t\not=0,\vspace{3pt}\\
H_0(\nabla_\xi H(\xi))=1, & \xi\in{\bf R}^N\setminus\{0\},\vspace{3pt}\\ 
H(\xi)\nabla_xH_0(\nabla_\xi H(\xi))=\xi, &  \xi\in{\bf R}^N,
\end{array}
\right.
\end{equation}
\begin{equation}
\label{eq:2.2}
\left\{
\begin{array}{ll}
x\cdot\nabla H_0(x)=H_0(x), & x\in{\bf R}^N,\vspace{3pt}\\
\nabla H_0(tx)=\mbox{sign}\,(t)\,\nabla H_0(x), & x\in{\bf R}^N\setminus\{0\},\,\,\,t\not=0,\vspace{3pt}\\
H(\nabla H_0(x))=1, & x\in{\bf R}^N\setminus\{0\},\vspace{3pt}\\
H_0(x)\nabla_\xi H(\nabla H_0(x))=x, & x\in{\bf R}^N.
\end{array}
\right.
\end{equation}
Here $\xi\cdot\nabla_\xi H(\xi)$, $H(\xi)\nabla_\xi H(\xi)$ 
and $x\cdot\nabla H_0(x)$, $H_0(x)\nabla H_0(x)$ are taken to be $0$ at $\xi=0$ and $x=0$, respectively.  
See \cite{CS} and \cite{FK}. 
\vspace{5pt}
\newline
{\bf Proof of Theorem~\ref{Theorem:1.1}.} 
Let $v(x)=v^\sharp(H_0(x))$, $v^*(x):=v^\sharp(|x|)$, $w(x)=w^\sharp(H_0(x))$ and $w^*(x):=w^\sharp(|x|)$
for $x\in{\bf R}^N$. Assume that $v^*$, $w^*\in C^2({\bf R}^N)$. 
It follows that 
\begin{equation}
\label{eq:2.3}
\begin{split}
 & (\nabla v)(x) =(\partial_rv^\sharp)(H_0(x))\nabla H_0(x),\\
 & H((\nabla v)(x))\nabla_\xi H((\nabla v)(x))
 =(\partial_r v^\sharp)(H_0(x))\nabla_\xi H(\nabla H_0(x))
=\frac{(\partial_rv^\sharp)(H_0(x))x}{H_0(x)},
\end{split}
\end{equation}
for $x\in{\bf R}^N\setminus\{0\}$. 
Then
\begin{equation}
\label{eq:2.4}
\begin{split}
 & (\Delta_H v)(x)=\mbox{div}\,\left(\frac{(\partial_rv^\sharp)(H_0(x))}{H_0(x)}x\right)\\
 & =\frac{[(\partial_r^2v^\sharp)(H_0(x))\nabla H_0(x)\cdot x+(\partial_rv^\sharp)(H_0(x))N]H_0(x)
 -(\partial_rv^\sharp)(H_0(x))x\nabla H_0(x)}{H_0(x)^2}\\
 & =\frac{(\partial_r^2v^\sharp)(H_0(x))H_0(x)^2+(N-1)(\partial_rv^\sharp)(H_0(x))H_0(x)}{H_0(x)^2}\\
 & =(\partial_r^2v^\sharp)(r)+\frac{N-1}{r}(\partial_rv^\sharp)(r)\,\biggr|_{r=H_0(x)}
 =(\Delta v^*)(y)\quad\mbox{with}\quad y=H_0(x)
\end{split}
\end{equation}
for $x\in{\bf R}^N\setminus\{0\}$. 
Since $(\partial_r v^\sharp)(0)=0$, 
it follows from \eqref{eq:2.3} that 
$$
H((\nabla v)(x))\nabla_{\xi_i}H((\nabla v)(x))=(\partial_r^2v^\sharp)(0)x_i+o(|x|)\quad\mbox{near $x=0$},
$$
where $i\in\{1,\dots,N\}$. 
This implies that 
\begin{equation}
\label{eq:2.5}
(\Delta_H v)(0)=N(\partial_r^2v^\sharp)(0)=(\Delta v^*)(0). 
\end{equation}
By \eqref{eq:2.4} and \eqref{eq:2.5} we obtain assertion~(a). 
Furthermore, we deduce from assertion~(a) that 
\begin{equation*}
\begin{split}
[\Delta_H(v+w)](x) & =[\Delta(v^*+w^*)](|x|)=(\Delta v^*)(|x|)+(\Delta w^*)(|x|)\\
 & =(\Delta_H v)(x)+(\Delta_H w)(x)
\end{split}
\end{equation*}
for $x\in{\bf R}^N$. Thus assertion~(b) follows, and the proof is complete.
$\Box$
\vspace{5pt}
\newline
Theorem~\ref{Theorem:1.3}  immediately follows from Theorem~\ref{Theorem:1.1}~(a) and (H3).
Similarly, 
we can find the following $H_0$-radially functions. 
\begin{itemize}
  \item 
  (Finsler Gauss kernel)\\
  Let
  $$
  G_{H_0}(x,t):=(4\pi t)^{-\frac{N}{2}}\exp\left(-\frac{H_0(x)^2}{4t}\right).
  $$
  Then $G_{H_0}$ is a solution of \eqref{eq:1.4} in ${\bf R}^N\times(0,\infty)$. 
  (See also \cite[Example~4.3]{OS1}.)
  \item 
  (Finsler Barenblatt solution)\\
  Let $m>1$ and $C>0$. Set 
  $$
  {\mathcal U}_{H_0}(x,t):=t^{-\alpha}(C-kH_0(x)^2 t^{-2\beta})_+^{\frac{1}{m-1}}, 
  $$
  where $(s)_+:=\max\{s,0\}$ and
  $$
  \alpha:=\frac{N}{N(m-1)+2},\qquad
  \beta:=\frac{\alpha}{N},\qquad
  k:=\frac{\alpha(m-1)}{2mN}.
  $$
  Then ${\mathcal U}_{H_0}$ is a solution of the Finsler porous medium equation 
  $$
  \partial_t v=\Delta_H v^m\quad\mbox{in}\quad{\bf R}^N\times(0,\infty).
  $$ 
    \item
  (Singular solutions to the $m$-th order Finsler-Laplace equation)\\
  Let $m\in\{1,2,\dots\}$. Set 
  $$
  v(x):=
  \left\{
  \begin{array}{ll}
  H_0(x)^{-N+2m} & \mbox{if}\quad N-2m\not\in\{-2h\,:\,h=0,1,2,\dots\},\vspace{3pt}\\
  H_0(x)^{-N+2m}\log H_0(x) & \mbox{if}\quad N-2m\in\{-2h\,:\,h=0,1,2,\dots\}.
  \end{array}
  \right.
  $$
  Then $v$ satisfies 
  $$
  (\Delta_H)^m v=c\delta
  $$
  for some constant $c\in{\bf R}$, where $\delta$ is the Dirac measure in ${\bf R}^N$. 
  \item 
  (Finsler Talenti function)\\
  Let $p>1$. Set 
  $$
  w(x):=\left(A+BH_0(x)^{\frac{p}{p-1}}\right)^{1-\frac{N}{p}},
  $$
  where $A>0$ and $B>0$.
  Then $w$ satisfies 
  $$
  -\Delta_H\,w=w^p\quad\mbox{in}\quad{\bf R}^N.
  $$
\end{itemize}
Due to assertion~(ii) of Theorem~\ref{Theorem:1.1}, 
we have a explicit representation of the $H_0$-radially symmetric solution of \eqref{eq:1.6}. 
Set
$$
I(z):=\int_{{\bf S}^{N-1}}e^{z\theta_1}\,d\theta. 
$$
\begin{theorem}
\label{Theorem:2.1}
Let $\varphi$ be an $H_0$-radially symmetric, bounded and continuous function in ${\bf R}^N$. 
Then 
\begin{equation}
\label{eq:2.6}
u(x,t):=(4\pi t)^{-\frac{N}{2}}\exp\left(-\frac{H_0(x)^2}{4t}\right)\int_0^\infty
I\left(\frac{H_0(x)r}{2t}\right)\exp\left(-\frac{r^2}{4t}\right)\varphi^\sharp(r)r^{N-1}\,dr
\end{equation}
is a solution of \eqref{eq:1.4} such that $u(x,0)=\varphi(x)$ in ${\bf R}^N$. 
\end{theorem}
{\bf Proof.}
Let $u^*(x,t):=e^{t\Delta}\varphi^*$, that is 
\begin{equation*}
\begin{split}
 u^*(x,t) &=(4\pi t)^{-\frac{N}{2}}\int_{{\bf R}^N} \exp\left(-\frac{|x-y|^2}{4t}\right)\varphi^*(y)\,dy\\
 & =(4\pi t)^{-\frac{N}{2}}\exp\left(-\frac{|x|^2}{4t}\right)
 \int_{{\bf R}^N}
 \exp \left( \dfrac{x \cdot y}{2t} \right)
 \exp \left( - \dfrac{|y|^2}{4t} \right)\varphi^\sharp(|y|) \, dy\\
 &= (4\pi t)^{-\frac{N}{2}}\exp\left(-\frac{|x|^2}{4t} \right)\\
 &\quad \times \int^\infty_0 \int_{{\bf S}^{N-1}}\exp \left(\frac{r|x|}{2t}\frac{x}{|x|}\theta\right)\,d\theta 
 \exp\left(-\frac{r^2}{4t}\right)\varphi^\sharp(r) r^{N-1}\,dr\\
 & =(4\pi t)^{-\frac{N}{2}}\exp\left(-\frac{|x|^2}{4t} \right)
 \int^\infty_0 I\left(\frac{r|x|}{2t}\right)\exp\left(-\frac{r^2}{4t}\right)\varphi^\sharp(r) r^{N-1}\,dr.
\end{split}
\end{equation*}
This together with Theorem~\ref{Theorem:1.1} implies that 
the function~$u$ defined by \eqref{eq:2.6} is a solution of \eqref{eq:1.4} such that $u(x,0)=\varphi(x)$ in ${\bf R}^N$. 
Thus Theorem~\ref{Theorem:2.1} follows.
$\Box$
\begin{remark}
\label{Remark:2.1}
Let $N=2$. Then
\begin{equation*}
\begin{split}
I(z) & =\int^{2\pi}_0 \exp(z\cos\theta)\,d\theta
=\frac{1}{i}\oint_{\Gamma} \exp \left(\frac{z}{2}(\xi + \xi^{-1})\right)\,\frac{d\xi}{\xi}\\
 & =2\pi\,\mathrm{Res}\left(\exp\left(\dfrac{z}{2}(\xi+\xi^{-1})\right) \, \xi^{-1} ; 0 \right),
 \quad z>0,
\end{split}
\end{equation*}
where $\Gamma$ denotes the unit circle centered at the origin in the complex plane. 
Noting that
\begin{align*}
 \exp\left(\frac{z}{2}(\xi+\xi^{-1})\right) \, \xi^{-1}
 &= \sum_{n=0}^\infty \sum_{m=0}^\infty \frac{1}{n!}\frac{1}{m!} 
 \left(\dfrac{z}{2}\right)^{n+m}\xi^{n-m-1}
\end{align*}
for $\xi\neq 0$, we see that
$$
\mathrm{Res}\left( \exp \left(\frac{z}{2}(\xi+\xi^{-1})\right)\,\xi^{-1} ; 0 \right)
= \sum_{n=0}^\infty \frac{1}{(n!)^2} \left(\frac{z}{2}\right)^{2n}
= I_0(z),
$$
where $I_0$ denotes the modified Bessel function of the first kind defined by
$$
I_0(z) := \sum_{n=0}^\infty \dfrac{(z/2)^{2n}}{(n!)^2} \quad \mbox{ for } \ z \in \mathbb{C} \setminus (-\infty,0].
$$
Therefore we deduce that $I(z)=2\pi I_0(z)$ for $z>0$. 
\end{remark}
\section{Cauchy-Dirichlet problem}
Let $R\geq1$, $\Omega=B_{H_0}(0,R)$ and $\phi\in C^\infty_0(\Omega)$. 
Consider the Cauchy-Dirichlet problem 
\begin{equation}
\label{eq:3.1}
\left\{
\begin{array}{ll}
\partial_t u=\Delta_H\, u & \mbox{in}\quad \Omega\times(0,\infty),\vspace{3pt}\\
u=0 & \mbox{on}\quad \partial\Omega\times(0,\infty),\vspace{3pt}\\
u(\cdot,0)=\phi & \mbox{in}\quad\Omega.
\end{array}
\right.
\end{equation}
Problem~\eqref{eq:3.1} possesses a unique solution $u$ in $\Omega\times[0,\infty)$ such that 
\begin{equation}\label{eq:3.2}
u\in W^{1,2}(0,\infty:L^2(\Omega)) \cap C([0,\infty): H_0^1(\Omega)). 
\end{equation}
To prove this fact, we define a functional $\psi$ on $\mathcal H := L^2(\Omega)$ by
\begin{equation}
\label{eq:3.3}
\psi(u) :=
\begin{cases}
\displaystyle{\frac{1}{2}\int_\Omega H(\nabla u(x))^2 \, d x} \ &\mbox{ if } \ u \in H^1_0(\Omega),\\
+\infty \ &\mbox{ otherwise.}
\end{cases}
\end{equation}
Then $\psi$ is convex and $D(\psi) := \{w \in \mathcal H \colon \psi(w) < +\infty\} = H^1_0(\Omega)$. 
Moreover, $\psi$ is lower semicontinuous in $\mathcal H$. 
Indeed, let $\lambda \geq 0$ be arbitrarily fixed and let $u \in \mathcal H$ and 
$u_n \in [\psi \leq \lambda] := \{w \in \mathcal H \colon \psi(u_n) \leq \lambda\}$ be such that 
$u_n \to u$ strongly in $\mathcal H$. 
Then $\{u_n\}$ is bounded in $H^1_0(\Omega)$ by the equivalence of $H(\cdot)$ to the Euclidean norm of $\mathbb{R}^N$. 
Hence $u_n \to u$ weakly in $H^1_0(\Omega)$. 
Now, since the set $[\psi \leq \lambda]$ is convex and closed in the strong topology of $H^1_0(\Omega)$, 
so is it in the weak topology. 
Hence $u$ belongs to $[\psi \leq \lambda]$. Thus $[\psi \leq \lambda]$ turns out to be closed in $\mathcal H$, and then,  
$\psi$ is lower semicontinuous in $\mathcal H$.
Next 
we define  the  subdifferential operator $\partial \psi : \mathcal H \to 2^{\mathcal H}$ of $\psi$ by
$$
\partial \psi(u) := \{\xi \in \mathcal H : \psi(v)-\psi(u) \geq (\xi,v-u)_{\mathcal H} \ \mbox{ for all } v \in  D(\psi)\},
$$
where $(\cdot,\cdot)_{\mathcal H}$ is  the  inner product of $\mathcal H$. 
Here we remark that $\partial \psi$ is maximal monotone in $\mathcal H$ 
since $\psi$ is proper (i.e., $\psi \not\equiv +\infty$), lower semicontinuous and convex in $\mathcal H$. 
We note by \eqref{eq:1.1} that
$\psi|_{H^1_0(\Omega)}$ is Fr\'echet differentiable and 
its Fr\'echet derivative $d\psi|_{H^1_0(\Omega)}(u)$ coincides with $-\Delta_H u$ (in $H^{-1}(\Omega)$) 
for $u \in H_0^1(\Omega)$. 
Moreover, by the definition of subdifferentials 
we see that 
$$
\partial \psi(u) \subset \partial \psi|_{H^1_0(\Omega)}(u) = \{d\psi|_{H^1_0(\Omega)}(u) \}
$$ 
for $u \in D(\partial \psi) := \{u \in H^1_0(\Omega) : \partial \psi(u) \neq \emptyset\}$. 
This implies that  
$$
\partial \psi(u) = \{- \Delta_H u\}\quad\mbox{in}\ \mathcal H
$$ 
for $u\in D(\partial \psi) = \{u \in H^1_0(\Omega) : \Delta_H u \in \mathcal H\}$. 
Hence \eqref{eq:3.3} is reduced to an abstract Cauchy problem in the Hilbert space $\mathcal H$ of the form
\begin{equation}\label{eq:3.4}
\partial_t u + \partial \psi(u) \ni 0 \ \mbox{ in } \mathcal H, \quad 0 < t < \infty, \quad u(0) = \phi,
\end{equation}
which has been well studied. In particular,  the 
global-in-time well-posedness for \eqref{eq:3.4} is guaranteed by K\=omura-Br\'ezis theory 
(see e.g.,~\cite[Theorem 3.6]{Br}). 
Thus we conclude that problem~\eqref{eq:3.1} has a unique strong solution $u$ satisfying \eqref{eq:3.2}.
\vspace{3pt}

In the rest of this section 
we improve the arguments in \cite[Section~4]{I} to prove the following proposition. 
\begin{proposition}
\label{Proposition:3.1}
Assume \eqref{eq:1.1} and \eqref{eq:1.3}. 
Let $u$ be a solution of \eqref{eq:3.1}, $R\geq1$ and $0<\ell<1/2$. 
Then there exist constants $C_*>0$ and $T_*\in(0,1)$ such that 
\begin{equation}
\label{eq:3.5}
\sup_{x\in\Omega}\int_{\Omega\,\cap \,B_{H_0}(x,1)}e^{-h(y,t)}|u(y,t)|\,dy
\le C_*\sup_{x\in\Omega}\int_{\Omega\,\cap \,B_{H_0}(x,1)}e^{-H_0(y)^2}|\phi(y)|\,dy
\end{equation}
for all $t\in(0,T_*]$, where $h(x,t):=H_0(x)^2(1+t^\ell)$. 
Here $C_*$ and $T_*$ depend only on $N$ and $\ell$. 
\end{proposition}
In what follows, the letter $C$ denotes generic positive constants independent of $\Omega$ 
and it may have different values also within the same line. 
We start with proving an energy estimate of the solution of \eqref{eq:3.1}. 
\begin{lemma}
\label{Lemma:3.1}
Assume the same conditions as in Proposition~{\rm\ref{Proposition:3.1}}. 
Set 
$$
z(x,t):=e^{-h(x,t)}u(x,t)=e^{-H_0(x)^2(1+t^\ell)}u(x,t).
$$
Then there exist $C>0$ and $T_1\in(0,1)$ such that 
\begin{equation}
\begin{split}
\label{eq:3.6}
 & \sup_{ t_1 < s < t}\int_{\Omega\,\cap\,B_{H_0}(x,R_2)} z(y,s)^2\zeta_x^2\,dy
+\iint_{Q_2}H(\nabla (z\zeta_x))^2\,dy\,ds\\
 & \qquad\qquad
\le C[(R_2-R_1)^{-2}+(t_1-t_2)^{-1}]\iint_{Q_2} z(y,s)^2\,dy\,ds
\end{split}
\end{equation}
and
\begin{equation}\label{eq:3.7}
\iint_{Q_1}z^{2\kappa}\,dy\,ds
\le C\left(\left[(R_2-R_1)^{-2}+(t_1-t_2)^{-1}\right]\iint_{Q_2}z^2\,dy\,ds\right)^\kappa, 
\end{equation}
for all $x\in\Omega$, $0<R_1<R_2$ and $0<t_2<t_1<t\le T_1$, where $\kappa:=(N+2)/N$ and 
\begin{equation}
\label{eq:3.8}
\begin{split}
 & \zeta_x(y,s):=\frac{1}{R_2-R_1}\min\{\max\{R_2-H_0(y-x),0\},R_2-R_1\}\\
 & \qquad\qquad\qquad\qquad
\times\frac{1}{t_1-t_2}\min\{\max\{s-t_2,0\},t_1-t_2\}
\end{split}
\end{equation}
for $x$, $y\in{\bf R}^N$ and $s>0$. 
Here 
$$
Q_1:=\left[\Omega\,\cap \,B_{H_0}(x,R_1)\right]\times(t_1,t),\qquad
Q_2:=\left[\Omega\,\cap \,B_{H_0}(x,R_2)\right]\times(t_2,t). 
$$  
\end{lemma}
{\bf Proof.}
It follows from \eqref{eq:2.2} and \eqref{eq:3.8} that 
\begin{equation}
\label{eq:3.9}
\begin{split}
 & 0\le\zeta_x\le 1\quad \mbox{in} \quad {\bf R}^N\times(0,\infty),
\qquad \zeta_x=1\quad\mbox{on}\quad Q_1,
\qquad \zeta_x=0\quad\mbox{on}\quad \partial_p Q_2,\vspace{3pt}\\
 & H(\nabla\zeta_x)\le\frac{1}{R_2-R_1}
\quad\mbox{and}\quad
0\le\partial_t\zeta_x\le\frac{1}{t_1-t_2}
\quad\mbox{in}\quad Q_2.
\end{split}
\end{equation}
Thanks to \eqref{eq:1.1}, \eqref{eq:1.2} and \eqref{eq:2.1}, 
by \eqref{eq:3.9} we have 
\begin{equation}
\label{eq:3.10}
\begin{split}
 & \iint_{Q_2}H(\nabla u)\nabla_\xi H(\nabla u)\nabla[e^{-2h}u\zeta_x^2]\,dy\,ds\\
 & =\iint_{Q_2}H(\nabla u)\nabla_\xi H(\nabla u)[-2e^{-2h}\nabla h u\zeta_x^2+e^{-2h}\nabla u\zeta_x^2+2e^{-2h}u\zeta_x\nabla\zeta_x]\,dy\,ds\\
 & \ge\iint_{Q_2}e^{-2h}H(\nabla u)^2\zeta_x^2\,dy\,ds
 -2\iint_{Q_2}e^{-2h}H(\nabla u)H_0(\nabla_\xi H(\nabla u))H(\nabla h)u\zeta_x^2\,dy\,ds\\
 & \qquad\qquad\qquad\qquad\qquad\,\,
  -2\iint_{Q_2}e^{-2h}H(\nabla u)H_0(\nabla_\xi H(\nabla u))H(\nabla\zeta_x)u\zeta_x\,dy\,ds\\
 & =\iint_{Q_2}e^{-2h}H(\nabla u)^2\zeta_x^2\,dy\,ds
 -2\iint_{Q_2}e^{-2h}H(\nabla u)H(\nabla h)u\zeta_x^2\,dy\,ds\\
 & \qquad\qquad\qquad\qquad\qquad\,\,
  -2\iint_{Q_2}e^{-2h}H(\nabla u)H(\nabla\zeta_x)u\zeta_x\,dy\,ds\\
 & \ge\frac{1}{2}\iint_{Q_2}e^{-2h}H(\nabla u)^2\zeta_x^2\,dy\,ds
 -4\iint_{Q_2}e^{-2h}[H(\nabla h)^2\zeta_x^2+H(\nabla\zeta_x)^2]u^2\,dy\,ds.
\end{split}
\end{equation}
Since $H$ is a norm of  ${\bf R}^N$, 
it follows that 
\begin{equation*}
\begin{split}
H(\nabla(z\zeta_x))^2
 & \le [e^{-h}u\zeta_x H(\nabla h)+e^{-h}H(\nabla u)\zeta_x+e^{-h}H(\nabla\zeta_x)u]^2\\
 & \le 3e^{-2h}[u^2\zeta_x^2H(\nabla h)^2+H(\nabla u)^2\zeta_x^2+H(\nabla\zeta_x)^2u^2].
\end{split}
\end{equation*}
These imply that 
\begin{equation*}
\begin{split}
 & \iint_{Q_2}H(\nabla u)\nabla_\xi H(\nabla u)\nabla[e^{-2h}u\zeta_x^2]\,dy\,ds\\
 & \ge\frac{1}{6}\iint_{Q_2}H(\nabla(z\zeta_x))^2\,dy\,ds
 -C\iint_{Q_2}[H(\nabla h)^2\zeta_x^2+H(\nabla\zeta_x)^2]z^2\,dy\,ds.
\end{split}
\end{equation*}
Therefore, 
multiplying \eqref{eq:1.4} by $e^{-2h}u\zeta_x^2$  and integrating it on $Q_2$, 
we obtain
\begin{equation}
\label{eq:3.11}
\begin{split}
 & \frac{1}{2}\int_{\Omega\,\cap\,B_{H_0}(x,R_2)} z(y,t)^2\zeta_x(y,t)^2\,dy
+\iint_{Q_2}  z^2[(\partial_t  h ) \zeta_x^2-\zeta_x\partial_t\zeta_x]\,dy\,ds\\
 & +\frac{1}{6}\iint_{Q_2}H(\nabla (z\zeta_x))^2\,dy\,ds
\le C\iint_{Q_2}[H(\nabla h)^2\zeta_x^2+H(\nabla\zeta_x)^2]z^2\,dy\,ds.
\end{split}
\end{equation}

Let $0<T_1\le 1$. Since $h(x,t):=H_0(x)^2(1+t^\ell)$, 
it follows from \eqref{eq:2.2} that 
\begin{equation}
\label{eq:3.12}
\partial_t h\ge\ell T_1^{\ell-1}H_0(x)^2,
\qquad
H(\nabla h)^2
=4H_0(x)^2(1+t^\ell)^2
\le 16H_0(x)^2,
\end{equation}
for all $x\in{\bf R}^N$ and $0<t<T_1$. 
Since $0<\ell<1/2$, 
by \eqref{eq:3.11} and \eqref{eq:3.12}, 
taking a sufficiently small $T_1>0$ if necessary, 
we obtain 
\begin{equation*}
\begin{split}
 & \frac{1}{2}\int_{\Omega\,\cap\,B_{H_0}(x,R_2)} z(y,t)^2\zeta_x^2\,dy
+\frac{1}{6}\iint_{Q_2}H(\nabla (z\zeta_x))^2\,dy\,ds\\
 & \le \iint_{Q_2}[ C  H(\nabla\zeta_x)^2+\zeta_x \partial_t\zeta_x]z^2\,dy\,ds,
\end{split}
\end{equation*}
which together with \eqref{eq:3.9} implies \eqref{eq:3.6}. 

Since $H$ is an equivalent norm to the Euclidean norm of ${\bf R}^N$,
by the Gagliardo-Nirenberg inequality we have 
\begin{equation*}
\begin{split}
 & \int^t_{ t_1}\int_{\Omega\,\cap\, B_{H_0}(x,R_2)}(z\zeta_x)^{2\kappa}\,dy\,ds\\
 & \le C\int^t_{ t_1}\int_{\Omega\,\cap\, B_{H_0}(x,R_2)}|\nabla(z\zeta_x)|^2\,dy\left(\int_{ \Omega \cap  B_{H_0}(x,R_2)}(z\zeta_x)^2\,dy\right)^{\frac{2}{N}}\,ds\\
 & \le C\int^t_{ t_1}\int_{\Omega\,\cap\, B_{H_0}(x,R_2)} H(\nabla(z\zeta_x))^2\,dy\,ds  \sup_{t_1 < s < t}\left(\int_{\Omega \cap B_{H_0}(x,R_2)}(z\zeta_x)^2(y,s)\,dy\right)^{\frac{2}{N}}
\end{split}
\end{equation*}
for all $x\in\Omega$ and $ t > t_1$. 
This together with  \eqref{eq:3.6} and \eqref{eq:3.9}  implies the desired inequality 
and Lemma~\ref{Lemma:3.1} follows.  
$\Box$\vspace{5pt}

By Lemma~\ref{Lemma:3.1} and the Besicovitch covering theorem 
we prove the following lemma.
\begin{lemma}
\label{Lemma:3.2}
Assume the same conditions as in Lemma~{\rm\ref{Lemma:3.1}}. 
Then, for any $\delta\in(0,1]$, 
there exists a positive constant $C$ such that 
\begin{equation}
\label{eq:3.13}
\sup_{x\in\Omega}\|z(t)\|_{L^{1+\delta}(\Omega\,\cap\,B_{H_0}(x,1))}
\le Ct^{-\frac{N}{2}\left(1-\frac{1}{1+\delta}\right)}\sup_{0<s<t}\sup_{x\in\Omega}\|z(s)\|_{L^1(\Omega\,\cap\,B_{H_0}(x,1))}
\end{equation}
for all $0<t<T_1$, where 
$T_1$ is as in Lemma~{\rm\ref{Lemma:3.1}}.
\end{lemma}
{\bf Proof.}
Let $x\in\Omega$, $\rho>0$ and $0<t<T_1$. 
Set 
$$
\sigma_n:=\sum_{i=1}^n 2^{-i},
\qquad
Q_n:=\left[\Omega\,\cap\,B_{H_0}(x,(1+\sigma_n)\rho)\right]\times\left(\frac{t}{4}\left(1-\frac{1}{2}\sigma_n\right),t\right). 
$$
It follows that 
$$
Q_0:=[\Omega\,\cap\,B_{H_0}(x,\rho)]\times\left( \frac{t}{4},t\right)
\subset Q_n\subset Q_{n+1}\subset Q_{\infty}=[\Omega\,\cap\, B_{H_0}(x,2\rho)]\times\left(\frac{t}{8} ,t\right)
$$
for $n=1,2,\dots$. 
Then we apply Lemma~\ref{Lemma:3.1} with $Q_1$ and $Q_2$ replaced by $Q_n$ and $Q_{n+1}$, respectively, 
and obtain 
\begin{equation}
\label{eq:3.14}
\|z\|_{L^{2\kappa}(Q_n)}\le CD_n^{1/2}\|z\|_{L^2(Q_{n+1})},
\end{equation}
where $D_n:=2^{2n}(\rho^{-2}+t^{-1})$. 
Furthermore, for any $\epsilon>0$, 
the H\"older  and Young inequalities imply  that 
\begin{equation}
\label{eq:3.15}
\begin{split}
CD_n^{1/2}\|z\|_{L^2(Q_{n+1})}
 & \le CD_n^{1/2}\|z\|_{L^1(Q_{n+1})}^\theta \|z\|_{L^{2\kappa}(Q_{n+1})}^{1-\theta}\\
 & \le\epsilon\|z\|_{L^{2\kappa}(Q_{n+1})}
 +\epsilon^{-(1-\theta)/\theta}\left[CD_n^{1/2}\|z\|_{L^1(Q_{n+1})}^{\theta}\right]^{\frac{1}{\theta}}\\
 & \le\epsilon\|z\|_{L^{2\kappa}(Q_{n+1})}
 +\epsilon^{-(1-\theta)/\theta}C^{1/\theta}D_n^{1/2\theta}\|z\|_{L^1(Q_{n+1})},
\end{split}
\end{equation}
where 
$$
\frac{1}{2}=\theta+\frac{1-\theta}{2\kappa}.
$$
By \eqref{eq:3.14} and \eqref{eq:3.15} we see that 
$$
\|z\|_{L^{2\kappa}(Q_n)}
\le \epsilon\|z\|_{L^{2\kappa}(Q_{n+1})}
+C\epsilon^{-(1-\theta)/\theta}D_n^{1/2\theta}\|z\|_{L^1(Q_{n+1})}
$$
for $n=0,1,2,\dots$. 
Then it follows that
\begin{equation}
\label{eq:3.16}
\begin{split}
\|z\|_{L^{2\kappa}(Q_0)}
 & \le\epsilon^k\|z\|_{L^{2\kappa}(Q_k)}
+C\epsilon^{-(1-\theta)/\theta}\sum_{n=0}^{k-1}\epsilon^nD_n^{1/2\theta}\|z\|_{L^1(Q_{n+1})}\\
 & \le\epsilon^k\|z\|_{L^{2\kappa}(Q_k)}
+C\epsilon^{-(1-\theta)/\theta}\sum_{n=0}^{k-1}\epsilon^n\left[2^{2n}(\rho^{-2}+t^{-1})\right]^{1/2\theta}\|z\|_{L^1(Q_{n+1})}
\end{split}
\end{equation}
for $k=1,2,\dots$. 
Taking a sufficiently small $\epsilon>0$ so that $\epsilon 2^{1/\theta}\le 1/2$ 
and passing to the limit, 
by \eqref{eq:3.16} we obtain
\begin{equation}
\label{eq:3.17}
\begin{split}
 & \|z\|_{L^{2\kappa}([\Omega\,\cap\, B_{H_0}(x,\rho)]\times(\frac{t}{4},t))}=\|z\|_{L^{2\kappa}(Q_0)}\\
 & \le C[\rho^{-2}+t^{-1}]^{1/2\theta}\|z\|_{L^1(Q_\infty)}
=C[\rho^{-2}+t^{-1}]^{ 1/2\theta}\|z\|_{L^1(B_{H_0}(x,2\rho)\times(\frac{t}{8},t))}
\end{split}
\end{equation}
for $x\in\Omega$, $\rho>0$ and $0<t<T_1$. 

On the other hand, 
by the H\"older inequality and \eqref{eq:3.6} we have
\begin{equation}
\label{eq:3.18}
\begin{split}
 & \|z(t)\|_{L^{1+\delta}(\Omega\,\cap\, B_{H_0}(x,\sqrt{t}))}
\le|B_{H_0}(x,\sqrt{t})|^{\frac{1}{1+\delta}-\frac{1}{2}}\|z(t)\|_{L^2(\Omega\,\cap\, B_{H_0}(x,\sqrt{t}))}\\
 & \le Ct^{\frac{N}{2}\left(\frac{1}{1+\delta}-\frac{1}{2}\right)}t^{-\frac{1}{2}}
 \left(\int_{t/4}^t\int_{\Omega\,\cap\, B_{H_0}(x,2\sqrt{t})}z^2\,dy\,ds\right)^{\frac{1}{2}}\\
 & \le Ct^{\frac{N}{2}\left(\frac{1}{1+\delta}-\frac{1}{2}\right)}t^{-\frac{1}{2}}
 \left(t|B_{H_0}(x,\sqrt{t})|\right)^{\frac{1}{2}(1-\frac{1}{\kappa})}
 \left(\int_{t/4}^t\int_{\Omega\,\cap\, B_{H_0}(x,2\sqrt{t})}z^{2\kappa}\,dy\,ds\right)^{\frac{1}{2\kappa}}.
 \end{split}
\end{equation}
Applying \eqref{eq:3.17} with $\rho=2\sqrt{t}$ to \eqref{eq:3.18}, 
we obtain 
\begin{equation}
\label{eq:3.19}
\begin{split}
 & \|z(t)\|_{L^{1+\delta}(\Omega\,\cap\, B_{H_0}(x,\sqrt{t}))}\\
 & \le Ct^{\frac{N}{2}\left(\frac{1}{1+\delta}-\frac{1}{2}\right)}t^{-\frac{1}{2}}
(t^{1+\frac{N}{2}})^{\frac{1}{2}(1-\frac{1}{\kappa})}t^{-\frac{1}{2\theta}}\int_{t/8}^t\int_{\Omega\,\cap\, B_{H_0}(x,4\sqrt{t})}|z|\,dy\,ds\\
 & \le Ct^{-\frac{N}{2}\left(1-\frac{1}{1+\delta}\right)}
\sup_{t/8<s<t}
\int_{\Omega\,\cap\,B_{H_0}(x,4\sqrt{t})}|z(y,s)|\,dy
 \end{split}
\end{equation}
for all $x\in\Omega$ and $0<t<T_1$. 

Let $x_0\in\Omega$. 
By the Besicovitch covering theorem (see e.g., \cite{EG}) 
we can find 
\begin{equation*}
\begin{split}
 & \mathcal{G}_1:=\{\overline{B_{H_0}(y_{1,k},\sqrt{t})}\},\quad \dots,\quad \mathcal{G}_n:=\{\overline{B_{H_0}(y_{n,k},\sqrt{t})}\}\\
 & \qquad\qquad\qquad\qquad\qquad
\subset{\mathcal F}:=\left\{\overline{B_{H_0}(y,\sqrt{t})}\,:\,y\in\overline{B_{H_0}(x_0,1)}\right\}
\end{split}
\end{equation*}
such that each $\mathcal{G}_i$ $(i=1,\dots,n)$ is a countable correction of disjoint balls in ${\mathcal F}$ 
and
\begin{equation}
\label{eq:3.20}
\overline{B_{H_0}(x_0,1)}\subset\bigcup_{i=1}^n\bigcup_k \overline{B_{H_0}(y_{i,k},\sqrt{t})},
\end{equation}
where $n$ is an integer depending only on $N$. 
Then there exists an integer $m$ depending only on $N$ such that 
\begin{equation}
\label{eq:3.21}
{}^\#\left\{y_{j,\ell}\,:\, B_{H_0}(y_{j,\ell},4\sqrt{t})\,\cap\,B_{H_0}(y_{i,k},4\sqrt{t})\not=\emptyset\right\}\le m
\end{equation}
for any $y_{i,k}$. 
Therefore, by \eqref{eq:3.19}, \eqref{eq:3.20} and \eqref{eq:3.21} we obtain  
\begin{equation*}
\begin{split}
\|z(t)\|_{L^{1+\delta}(\Omega\,\cap \,B_{H_0}(x_0,1))}
 & \le\left\|\sum_{i=1}^n\sum_k |z(t)|\chi_{ B_{H_0}(y_{i,k},\sqrt{t})}\right\|_{L^{1+\delta}(\Omega\,\cap \,B_{H_0}(x_0,1))}\\
 & \le\sum_{i=1}^n\sum_k\|z(t)\|_{L^{1+\delta}(\Omega\,\cap \,B_{H_0}(y_{i,k},\sqrt{t}))}\\
 & \le Ct^{-\frac{N}{2}\left(1-\frac{1}{1+\delta}\right)}\sum_{i=1}^n\sum_k
 \sup_{t/8<s<t}\int_{\Omega\,\cap \,B_{H_0}(y_{i,k},4\sqrt{t})}|z(y,s)|\,dy\\
 & \le Ct^{-\frac{N}{2}\left(1-\frac{1}{1+\delta}\right)}
 \sup_{t/8<s<t}\int_{\Omega\,\cap \,B_{H_0}(x_0,1+4\sqrt{t})}|z(y,s)|\,dy. 
\end{split}
\end{equation*}
Since $x_0$ is arbitrary, 
we deduce that 
\begin{equation*}
\begin{split}
\sup_{x\in\Omega}\|z(t)\|_{L^{1+\delta}(\Omega\,\cap \,B_{H_0}(x,1))}
 & \le Ct^{-\frac{N}{2}\left(1-\frac{1}{1+\delta}\right)}
\sup_{0<s<t}\sup_{x\in\Omega}\int_{\Omega\,\cap \,B_{H_0}(x,1+4\sqrt{t})}|z(y,s)|\,dy\\
 & \le Ct^{-\frac{N}{2}\left(1-\frac{1}{1+\delta}\right)}
\sup_{0<s<t}\sup_{x\in\Omega}\int_{\Omega\,\cap \,B_{H_0}(x,1)}|z(y,s)|\,dy
\end{split}
\end{equation*}
for all $0<t<T_1$. Thus \eqref{eq:3.13} holds. The proof is complete. 
$\Box$
\vspace{5pt}

Next, applying a similar argument as in \cite[Lemma~4.2]{I}, we have: 
\begin{lemma}
\label{Lemma:3.3}
Assume the same conditions as in Lemma~{\rm\ref{Lemma:3.1}}. 
Let $x\in\Omega$ and $\eta(y):=\min\{\max\{2-H_0(y-x),0\},1\}$. 
For $\epsilon>0$, $\sigma>0$ and $\delta>0$, 
set 
\begin{equation*}
\begin{split}
 & J_\pm^\epsilon(x,t):=\int_0^t\int_{\Omega\,\cap\,B_{H_0}(x,2)}
 s^\sigma e^{-h}\frac{H(\nabla u_\pm)^2}{u_\pm+\epsilon}[e^{-h}(u_\pm+\epsilon)]^\delta\eta^2\,dy\,ds,\\
 & u_\pm:=\max\{\pm u,0\},
 \qquad
 z_\pm=e^{-h}u_\pm.
\end{split}
\end{equation*}
Then, for any $\sigma\in(0,1/2)$ and $\delta\in(0,2\sigma/N)$,  
there exist positive constants $C$ and $T_2$ such that 
\begin{equation}
\label{eq:3.22}
\limsup_{\epsilon\to 0}J_\pm^\epsilon(x,t)\le C\int_0^t\int_{\Omega\,\cap\,B_{H_0}(x,2)} s^{\sigma-1}z_\pm^{1+\delta}\,dy\,ds
\end{equation}
for all $0<t<T_2$ and $x\in\Omega$.
\end{lemma}
{\bf Proof.}
For any $\epsilon>0$, set
$z_\pm^\epsilon := e^{-h}(u_\pm+\epsilon)$ and $\psi_\pm^\epsilon:=s^\sigma e^{-h} (z_\pm^\epsilon)^\delta\eta^2$. 
Then 
\begin{equation}
\label{eq:3.23}
\begin{split}
 & \lim_{\epsilon\to 0}\int_0^t\int_{\Omega\,\cap \,B_{H_0}(x,2)}(\partial_t u)\psi_\pm^{\epsilon}\,dy\,ds\\
 & =\lim_{\epsilon\to 0}\int_0^t\int_{\Omega\,\cap \,B_{H_0}(x,2)}[e^h\partial_t z_\pm+(u_\pm+\epsilon)\partial_t h]s^\sigma e^{-h}(z_\pm^\epsilon)^\delta\eta^2\,dy\,ds\\
 & =\int_0^t\int_{\Omega\,\cap \,B_{H_0}(x,2)}[s^{\sigma} (\partial_t z_\pm) z_\pm^\delta \eta^2 + s^\sigma z_\pm^{1+\delta} \eta^2 \partial_t h]\,dy\,ds\\
 & =\frac{1}{1+\delta}\int_0^t\int_{\Omega\,\cap \,B_{H_0}(x,2)}s^\sigma\eta^2\partial_t(z_\pm^{1+\delta})\,dy\,ds
	+\int_0^t\int_{\Omega\,\cap \,B_{H_0}(x,2)}s^\sigma z_\pm^{1+\delta}\eta^2\partial_t h\,dy\,ds\\
 & \ge -\frac{\sigma}{1+\delta}\int_0^t\int_{\Omega\,\cap \,B_{H_0}(x,2)} s^{\sigma-1}z_\pm^{1+\delta}\eta^2\,dy\,ds
+\int_0^t\int_{\Omega\,\cap \,B_{H_0}(x,2)}s^\sigma z_\pm^{1+\delta}\eta^2 \partial_t h\,dy\,ds.
\end{split}
\end{equation}
Furthermore, 
by \eqref{eq:1.2}, \eqref{eq:2.1} and \eqref{eq:2.2} 
we obtain 
\begin{equation}
\label{eq:3.24}
\begin{split}
 & \int_0^t\int_{\Omega\,\cap \,B_{H_0}(x,2)}H(\nabla u_\pm)\nabla_\xi H(\nabla u_\pm)\nabla\psi_\pm^\epsilon\,dy\,ds \\
 & =\int_0^t\int_{\Omega\,\cap \,B_{H_0}(x,2)}H(\nabla u_\pm)\nabla_\xi H(\nabla u_\pm)\nabla(s^\sigma e^{-h} (z_\pm^{\epsilon})^{\delta} \eta^2)\,dy\,ds \\
 & \ge \frac{\delta}{2} \int_0^t\int_{\Omega\,\cap \,B_{H_0}(x,2)}s^\sigma e^{-2h} H(\nabla u_\pm)^2 (z_\pm^\epsilon)^{-1+\delta} \eta^2 \,dy\,ds\\
 & \qquad\qquad 
 -C\int_0^t\int_{\Omega\,\cap \,B_{H_0}(x,2)}s^\sigma (z_\pm^\epsilon)^{1+\delta}[H(\nabla h)^2\eta^2+H(\nabla \eta)^2]\,dy\,ds\\
 & \ge \frac{\delta}{2} \int_0^t\int_{\Omega\,\cap \,B_{H_0}(x,2)}s^\sigma e^{-h}\frac{H(\nabla u_\pm)^2}{u_\pm+\epsilon}
 [e^{-h}(u_\pm+\epsilon)]^\delta\eta^2\,dy\,ds \\
 & \qquad \qquad 
 -C\int_0^t\int_{\Omega\,\cap \,B_{H_0}(x,2)}s^\sigma (z_\pm^\epsilon)^{1+\delta}[H(\nabla h)^2\eta^2+1]\,dy\,ds.
\end{split}
\end{equation}
On the other hand, by \cite[Chapter~II, Section~1]{DB02} 
we see that $u_+$ and $u_-$ are subsolutions of \eqref{eq:3.1}. 
Then, similarly to \eqref{eq:3.11}, 
we deduce from \eqref{eq:3.23} and \eqref{eq:3.24} that 
\begin{equation}
\label{eq:3.25}
\begin{split}
 & \frac{\delta}{2}\limsup_{\epsilon\to 0}\int_0^t\int_{\Omega\,\cap\,B_{H_0}(x,2)}
 s^\sigma e^{-h}\frac{H(\nabla u_\pm)^2}{u_\pm+\epsilon}[e^{-h}(u_\pm+\epsilon)]^\delta\eta^2\,dy\,ds\\
 & \le\frac{\sigma}{1+\delta}\int_0^t\int_{\Omega\,\cap\,B_{H_0}(x,2)} s^{\sigma-1}z_\pm^{1+\delta}\eta^2\,dy\,ds\\
 & \qquad\qquad
+\int_0^t\int_{\Omega\,\cap\,B_{H_0}(x,2)}s^\sigma z_\pm^{1+\delta}[CH(\nabla h)^2\eta^2-(\partial_t h)\eta^2+C]\,dy\,ds\\
 & \le C\int_0^t\int_{\Omega\,\cap\,B_{H_0}(x,2)} s^{\sigma-1}z_\pm^{1+\delta}\,dy\,ds\\
 & \qquad\qquad
 +\int_0^t\int_{\Omega\,\cap\,B_{H_0}(x,2)}s^\sigma z_\pm^{1+\delta}[CH(\nabla h)^2\eta^2-(\partial_t h)\eta^2]\,dy\,ds.
\end{split}
\end{equation}
Furthermore, by \eqref{eq:3.12} we can find $T_2\in(0,1)$ such that 
$$
CH(\nabla h)^2- \partial_t h\leq 0\quad\mbox{in}\quad{\bf R}^N\times(0,T_2).
$$
This together with \eqref{eq:3.25} implies \eqref{eq:3.22}. 
The proof is complete. 
$\Box$\vspace{5pt}

Now we are ready to prove Proposition~\ref{Proposition:3.1}.
\vspace{5pt}
\newline
{\bf Proof of Proposition~\ref{Proposition:3.1}.}
Let $x\in{\bf R}^N$ and let $\eta$ be as in Lemma~\ref{Lemma:3.3}. 
Let $\delta$ be a sufficiently small positive constant. 
Following \cite{I} and Lemma~\ref{Lemma:3.3}, we set
\begin{equation*}
\begin{split}
 & F(v):=\int_0^v \min\{\tau^\delta,1\}\,d\tau
=\left\{
\begin{array}{ll}
v^{1+\delta}/(1+\delta) & \mbox{for}\quad 0\le v\le 1,\vspace{3pt}\\
v-\delta/(1+\delta) & \mbox{for}\quad v>1, 
\end{array}
\right.\\
 & 
\varphi_\pm^\epsilon(y,s):=e^{-h}F'(e^{-h}(u_\pm+\epsilon))\eta^2,
\quad
z_\pm=e^{-h}u_\pm,
\quad
z_\epsilon=e^{-h}(u_\pm+\epsilon).
\end{split}
\end{equation*}
Then $F'(v)=v^\delta$ for $0\le v\le 1$ and $F'(v)=1$ for $v>1$. 
Taking a sufficiently small $\delta>0$ if necessary, 
we can assume that 
\begin{equation}
\label{eq:3.26}
F(v)\leq v\quad\mbox{for}\quad v \ge 0
\qquad\mbox{and}\qquad
v \leq 2F(v)\quad\mbox{for}\quad v \ge 1.
\end{equation}
Let $0<\sigma<1-\ell$. 
Then we have 
\begin{equation*}
\begin{split}
 & \lim_{\epsilon\to 0}\int_0^t\int_{\Omega\,\cap\,B_{H_0}(x,2)}
\partial_t u_\pm\varphi_\pm^\epsilon\,dy\,ds
 =\int_0^t\int_{\Omega\,\cap\,B_{H_0}(x,2)}[\partial_t z_\pm+z_\pm\partial_t h]F'(z_\pm)\eta^2\,dy\,ds\\
 & 
=\int_{\Omega\,\cap\,B_{H_0}(x,2)}F(z_\pm(s))\eta^2\,dy\biggr|_{s=0}^{s=t}\\
 & \qquad\quad
+\int_0^t\int_{\Omega\,\cap\,B_{H_0}(x,2)}\left[z_\pm^{1+\delta}\chi_{\{z_\pm\le 1\}}+z_\pm\chi_{\{z_\pm>1\}}\right](\partial_t h)\eta^2\,dy\,ds.
\end{split}
\end{equation*}
Furthermore, similarly to \eqref{eq:3.24}, we obtain 
\begin{equation*}
\begin{split}
 & \int_0^t\int_{\Omega\,\cap\,B_{H_0}(x,2)}
H(\nabla u_\pm)\nabla_\xi H(\nabla u_\pm)\nabla\varphi_\pm^\epsilon\,dy\,ds\\
 & =\int_0^t\int_{\Omega\,\cap\,B_{H_0}(x,2)}H(\nabla u_\pm)\nabla_\xi H(\nabla u_\pm)
 e^{-h}[-z_\pm^\epsilon \nabla h+e^{-h}\nabla u_\pm]\delta (z_\pm^\epsilon)^{-1+\delta}\eta^2\chi_{\{z_\pm^\epsilon\le 1\}}\,dy\,ds\\
 & \qquad
 +\int_0^t\int_{\Omega\,\cap\,B_{H_0}(x,2)}
 H(\nabla u_\pm)\nabla_\xi H(\nabla u_\pm)[-e^{-h}\nabla h \eta^2
 +e^{-h}2\eta\nabla\eta]F'(z_\pm^\epsilon)\,dy\,ds\\
 & \ge\delta\int_0^t\int_{\Omega\,\cap\,B_{H_0}(x,2)}
 e^{-2h}H(\nabla u_\pm)^2(z_\pm^\epsilon)^{-1+\delta}\eta^2\chi_{\{z_\pm^\epsilon\le 1\}}\,dy\,ds\\
 & \qquad 
 -\delta \int_0^t\int_{\Omega\,\cap\,B_{H_0}(x,2)}e^{-h}H(\nabla u_\pm)H(\nabla h)
 (z_\pm^\epsilon)^\delta\eta^2 \chi_{\{z_\pm^\epsilon\le 1\}}\,dy\,ds\\
 & \qquad
-\int_0^t\int_{\Omega\,\cap\,B_{H_0}(x,2)}e^{-h}H(\nabla u_\pm)(H(\nabla h)\eta^2+2\eta H(\nabla\eta))
[(z_\pm^\epsilon)^\delta\chi_{\{z_\pm^\epsilon\le 1\}}+\chi_{\{z_\pm^\epsilon>1\}}]\,dy\,ds\\
 & \ge -\int_0^t\int_{\Omega\,\cap\,B_{H_0}(x,2)}
 s^{\sigma}e^{-2h}H(\nabla u_\pm)^2 (z_\pm^\epsilon)^{-1+\delta}\eta^2\,dy\,ds\\
 & \qquad
 -C\int_0^t\int_{\Omega\,\cap\,B_{H_0}(x,2)}
 s^{-\sigma}(z_\pm^\epsilon)^{1+\delta}[H(\nabla h)^2\eta^2+H(\nabla\eta)^2]\chi_{\{z_\pm^\epsilon\le 1\}}\,dy\,ds\\
 & \qquad
 -C\int_0^t\int_{\Omega\,\cap\,B_{H_0}(x,2)}
 s^{-\sigma}(z_\pm^\epsilon)^{1-\delta}[H(\nabla h)^2\eta^2+H(\nabla\eta)^2]\chi_{\{z_\pm^\epsilon>1\}}\,dy\,ds.
\end{split}
\end{equation*}
Similarly to \eqref{eq:3.25}, 
these imply that 
\begin{equation}
\label{eq:3.27}
\begin{split}
 & \int_{\Omega\,\cap\,B_{H_0}(x,2)}F(z_\pm(s))\eta^2\,dy\biggr|_{s=0}^{s=t}\\
 & \le\limsup_{\epsilon\to 0}\int_0^t\int_{\Omega\,\cap\,B_{H_0}(x,2)}
s^{\sigma}e^{-2h}H(\nabla u_\pm)^2(z_\pm^\epsilon)^{-1+\delta}\eta^2\,dy\,ds\\
 & 
+C \int_0^t\int_{\Omega\,\cap\,B_{H_0}(x,2)}s^{-\sigma}[z_\pm^{1+\delta}\chi_{\{z_\pm\le 1\}}+z_\pm^{1-\delta}\chi_{\{z_\pm> 1\}}]
[H(\nabla h)^2\eta^2+H(\nabla\eta)^2]\,dy\,ds\\
 & 
 -\int_0^t\int_{\Omega\,\cap\,B_{H_0}(x,2)}[z_\pm^{1+\delta}\chi_{\{z_\pm\le 1\}}+z_\pm\chi_{\{z_\pm>1\}}](\partial_t h)\eta^2\,dy\,ds.
\end{split}
\end{equation}
On the other hand, since $\ell+\sigma<1$, we can find a constant $T_2\in(0,T_1)$ such that 
\begin{equation}
\label{eq:3.28}
\begin{split}
 & \partial_t h = \ell s^{\ell - 1}H_0(y)^2
 \ge\ell T_2^{\ell+\sigma-1}s^{-\sigma}H_0(y)^2,\\
 & Cs^{-\sigma} H(\nabla h)^2
 = 4C(1+s^\ell)^2 s^{-\sigma}H_0(y)^2
  \le 4C(1+T_2^\ell)^2 s^{-\sigma}H_0(y)^2 \le \partial_t h,
\end{split}
\end{equation}
for $y\in{\bf R}^N$ and $0<t<T_2$. 
Since $z_\pm^{1-\delta}\chi_{\{z_\pm>1\}}\le z_\pm\chi_{\{z_\pm>1\}}$, 
by \eqref{eq:3.27} and \eqref{eq:3.28} we obtain
\begin{equation*}
\begin{split}
 & \int_{\Omega\,\cap\,B_{H_0}(x,2)}F(z_\pm(s))\eta^2\,dy\biggr|_{s=0}^{s=t}\\
 & \le\limsup_{\epsilon\to 0}J_\pm^\epsilon(x,t)
+C\int_0^t\int_{\Omega\,\cap\,B_{H_0}(x,2)}s^{-\sigma}[z_\pm^{1+\delta}\chi_{\{z_\pm\le 1\}}+z_\pm\chi_{\{z_\pm>1\}}]H(\nabla\eta)^2\,dy\,ds\\
& \le\limsup_{\epsilon\to 0}J_\pm^\epsilon(x,t)
+C\int_0^t\int_{\Omega\,\cap\,B_{H_0}(x,2)}s^{-\sigma}F(z_\pm)\,dy\,ds
\end{split}
\end{equation*}
for $x\in\Omega$ and $0<t<T_2$. 
This together with Lemma~\ref{Lemma:3.3} implies that 
\begin{equation}
\label{eq:3.29}
\begin{split}
 & \int_{\Omega\,\cap\,B_{H_0}(x,1)}F(|z(t)|)\,dy
 \le\int_{\Omega\,\cap\,B_{H_0}(x,2)}F(|z(0)|)\,dy\\
 & \qquad
+C\int_0^t\int_{\Omega\,\cap\,B_{H_0}(x,2)} s^{\sigma-1}|z|^{1+\delta}\,dy\,ds
+C\int_0^t\int_{\Omega\,\cap\,B_{H_0}(x,2)} s^{-\sigma}F(z)\,dy\,ds
\end{split}
\end{equation}
for $x\in\Omega$ and $0<t<T_2$. 

Set 
$$
I(t):=\sup_{0\le s\le t}\sup_{x\in\Omega}\int_{\Omega\,\cap\,B_{H_0}(x,1)} |z(y,s)|\,dy,\qquad t>0. 
$$
If $I(0)=0$, then $u\equiv 0$ in $\Omega\times(0,\infty)$ and Proposition~\ref{Proposition:3.1} holds. 
So it suffices to consider the case $I(0)\not=0$. 
Furthermore, thanks to \eqref{eq:1.5}, 
we can assume, without loss of generality, that 
\begin{equation}
\label{eq:3.30}
I(0)=\sup_{x\in\Omega}\int_{\Omega\,\cap\,B_{H_0}(x,1)} |z(y,0)|\,dy=\sup_{x\in\Omega}\int_{\Omega\,\cap\,B_{H_0}(x,1)}e^{- H_0(y)^2}|\phi(y)|\,dy=1.
\end{equation}
On the other hand, 
for any $x\in{\bf R}^N$, there exist an integer $M$ depending only on $N$ 
and a set $\{x_j\}_{j=1}^M$ such that
\begin{equation}
\label{eq:3.31}
B_{H_0}(x,2)\subset \bigcup_{j=1}^MB_{H_0}(x_j,1)
\end{equation}
(see e.g., \cite[Lemma~2.1]{IS1}). 
Then, by \eqref{eq:3.26}, \eqref{eq:3.29}, \eqref{eq:3.30} and \eqref{eq:3.31}
we see that 
\begin{equation*}
\begin{split}
 & \sup_{x\in\Omega}\int_{\Omega\,\cap\,B_{H_0}(x,1)}F(|z(y,t)|)\,dy\\
 & \le\sup_{x\in\Omega}
 \int_{\Omega\,\cap\,B_{H_0}(x,2)}|z(y,0)|\,dy
 +C\sup_{x\in\Omega}\int_0^t\int_{\Omega\,\cap\,B_{H_0}(x,2)} s^{\sigma-1}|z(y,s)|^{1+\delta}\,dy\,ds\\
 & \qquad\qquad\qquad\qquad\qquad\qquad\,\,\,\,
 +C\sup_{x\in\Omega}\int_0^t\int_{\Omega\,\cap\,B_{H_0}(x,2)} s^{-\sigma}|z(y,s)|\,dy\,ds\\
 & \le M\sup_{x\in\Omega}
 \int_{\Omega\,\cap\,B_{H_0}(x,1)}|z(y,0)|\,dy
 +CM\sup_{x\in\Omega}\int_0^t\int_{\Omega\,\cap\,B_{H_0}(x,1)} s^{\sigma-1}|z(y,s)|^{1+\delta}\,dy\,ds\\
 & \qquad\qquad\qquad\qquad\qquad\qquad\,\,\,\,
 +CM\sup_{x\in\Omega}\int_0^t\int_{\Omega\,\cap\,B_{H_0}(x,1)} s^{-\sigma}|z(y,s)|\,dy\,ds
\end{split}
\end{equation*}
for $0<t<T_2$. 
This together with Lemma~\ref{Lemma:3.2} implies that
\begin{equation}
\label{eq:3.32}
\begin{split}
 & \sup_{x\in\Omega}\int_{\Omega\,\cap\,B_{H_0}(x,1)}F(|z(y,t)|)\,dy\\
 & \le MI(0)+C\int_0^t s^{\sigma-1}\sup_{x\in\Omega}\|z(s)\|_{L^{1+\delta}(\Omega\,\cap\,B_{H_0}(x,1))}^{1+\delta}\,ds
 +C\int_0^t s^{-\sigma}I(s)\,ds\\
 & \le M+C\int_0^t s^{\sigma-1-\frac{\delta N}{2}}I(s)^{1+\delta}\,ds+C\int_0^t s^{-\sigma}I(s)\,ds
\end{split}
\end{equation}
for all $0<t<T_2$. Since $I(0)=1<5M$, 
we can define 
$$
T_*:=\sup\left\{t\in(0,T_2]\,:\,I(t)\le L:=5M+|B_{H_0}(0,1)|\right\}. 
$$
 Since $\delta N/2 < \sigma < 1$ (see Lemma \ref{Lemma:3.3}),  taking a sufficiently small $T_2>0$ if necessary, 
by \eqref{eq:3.32} we obtain 
\begin{equation}
\label{eq:3.33}
\begin{split}
\sup_{x\in\Omega}\int_{\Omega\,\cap\,B_{H_0}(x,1)}F(|z( y, t)|)\,dy
\le M+CL^{1+\delta}\int_0^t s^{\sigma-1-\frac{\delta N}{2}}\,ds+CL\int_0^t s^{-\sigma}\,ds
\le 2M
\end{split}
\end{equation}
for all $0<t\le T_2$. 
On the other hand, 
it follows from \eqref{eq:3.26} that 
\begin{equation}
\label{eq:3.34}
\begin{split}
 & \sup_{x\in\Omega}\int_{\Omega\,\cap\,B_{H_0}(x,1)}|z(y,t)|\,dy\\
 & \le \sup_{x\in\Omega}\int_{\Omega\,\cap\,B_{H_0}(x,1)} \chi_{\{|z(y,t)|<1\}}(y)\,dy
 +\sup_{x\in\Omega}\int_{\Omega\,\cap\,B_{H_0}(x,1)} |z(y,t)|\chi_{\{|z(y,t)|\ge 1\}}(y)\,dy\\
 & \le |B_{H_0}(0,1)|+2\sup_{x\in{\bf R}^N}\int_{\Omega\,\cap\,B_{H_0}(x,1)} F(|z(y,t)|)\,dy
\end{split}
\end{equation}
for all $t>0$. 
Therefore we deduce from \eqref{eq:3.33} and \eqref{eq:3.34} that
$$
I(t)\le |B_{H_0}(0,1)|+4M
$$
for all $0<t\le T_*$. This means that $T_*=T_2$ and 
$I(t)\le L=LI(0)$
for all $0<t\le T_2$. 
Thus \eqref{eq:3.5} holds for $t\in(0,T_2]$ 
and Proposition~\ref{Proposition:3.1} follows.  
$\Box$ 
\vspace{5pt}

By Proposition~\ref{Proposition:3.1} we have: 
\begin{proposition}
\label{Proposition:3.2}
Assume the same conditions as in Proposition~{\rm\ref{Proposition:3.1}}. 
Then there exist positive constants $C_1$ and $\sigma'$ such that 
\begin{equation}
\label{eq:3.35}
\begin{split}
 & \sup_{x\in\Omega}\int_0^t\int_{\Omega\,\cap\,B_{H_0}(x,1)}e^{-h(y,s)}H(\nabla u(y,s))\,dy\,ds\\
 & \qquad\quad
 \le C_1 t^{\sigma'}\sup_{x\in\Omega}\int_{\Omega\,\cap\,B_{H_0}(x,1)}e^{-H_0(y)^2}|\phi(y)|\,dy
\end{split}
\end{equation}
for all $t\in(0,T_*)$, where $T_*$ is as in Proposition~{\rm\ref{Proposition:3.1}}. 
Furthermore, there exists a positive constant $C_2$ such that 
\begin{equation}
\label{eq:3.36}
\begin{split}
 & \sup_{x\in\Omega}\left(\int_{\Omega\,\cap\,B_{H_0}(x,1)}
|e^{-h(y,t)}u(y,t)|^2\,dy\right)^{\frac{1}{2}}\\
 & \qquad\quad
\le C_2t^{-\frac{N}{4}}
\sup_{x\in\Omega}\int_{\Omega\,\cap\,B_{H_0}(x,1)}e^{-H_0(y)^2}|\phi(y)|\,dy
\end{split}
\end{equation}
for all $t\in(0,T_*)$. 
\end{proposition}
{\bf Proof.}
We use the same notation  and assume $I(0)=1$  as in the proof of Proposition~\ref{Proposition:3.1}. 
By Lemma~\ref{Lemma:3.3} we have
\begin{equation*}
\begin{split}
 & \int_0^t\int_{\Omega\,\cap\,B_{H_0}(x,1)}e^{-h(y,s)}H(\nabla u(y,s))\,dy\,ds\\
 & \le\frac{1}{2}\limsup_{\epsilon\to 0}\,[J_+^\epsilon(x,t)+J_-^\epsilon(x,t)]
 +\frac{1}{2}\int_0^t\int_{\Omega\,\cap\,B_{H_0}(x,2)} s^{-\sigma}|z|^{1-\delta}\eta^2\,dy\,ds\\
 & \le C\int_0^t\int_{\Omega\,\cap\,B_{H_0}(x,2)} s^{\sigma-1}|z|^{1+\delta}\,dy\,ds
 +\frac{ M}{2}\int_0^t\int_{\Omega\,\cap\,B_{H_0}(x,1)} s^{-\sigma}[1+|z|]\,dy\,ds
\end{split}
\end{equation*}
for all $x\in\Omega$ and $0<t<T_*=T_2$. 
Then, similarly to \eqref{eq:3.32}, 
by Lemma~\ref{Lemma:3.2} and \eqref{eq:3.31} we obtain 
\begin{equation*}
\begin{split}
 & \sup_{x\in\Omega}\int_0^t\int_{\Omega\,\cap\,B_{H_0}(x,1)}e^{-h(y,s)}H(\nabla u(y,s))\,dy\,ds\\
 & \le C\int_0^t s^{\sigma-1-\frac{\delta N}{2}}I(s)^{1+\delta}\,ds
+\frac{ M}{2}\int_0^t s^{-\sigma}[I(s)+ |B_{H_0}(x,1)| ]\,ds
\end{split}
\end{equation*}
for all $0<t<T_*=T_2$. 
This together with Proposition~\ref{Proposition:3.1} implies that 
$$
\sup_{x\in\Omega}\int_0^t\int_{\Omega\,\cap\,B_{H_0}(x,1)}e^{-h(y,s)}H(\nabla u(y,s))\,dy\,ds
\le C\left[t^{\sigma-\frac{\delta N}{2}}+t^{1-\sigma}+t\right]I(0)
\le Ct^{\sigma'}I(0)
$$
for all $0<t<T_*=T_2$, where $\sigma'=\min\{\sigma-\delta N/2,1-\sigma,1\}$. 
Thus \eqref{eq:3.35} holds. 

On the other hand, 
by Lemma~\ref{Lemma:3.2} and Proposition~\ref{Proposition:3.1} 
we have
$$
\sup_{x\in\Omega}\|z(t)\|_{L^2(\Omega\,\cap\,B_{H_0}(x,1))}
\le Ct^{-\frac{N}{4}}\sup_{0<s<t}\sup_{x\in\Omega}\|z(s)\|_{L^1(\Omega\,\cap\,B_{H_0}(x,1))}
\le Ct^{-\frac{N}{4}}I(0)
$$
for all $0<t<T_*$. This implies \eqref{eq:3.36}. 
Thus Proposition~\ref{Proposition:3.2} follows. 
$\Box$
\section{Proof of Theorem~\ref{Theorem:1.2}}
In this section we prove Theorem~\ref{Theorem:1.2}. 
We prepare the following proposition.
\begin{proposition}
\label{Proposition:4.1}
Assume \eqref{eq:1.1} and \eqref{eq:1.3}. 
Let $u$ be a solution of \eqref{eq:3.1}. 
Then
\begin{equation}
\label{eq:4.1}
\sup_{0<t<S_\lambda}\int_{\Omega}e^{-2g(y, t )}u(y, t )^2\,dy
\le \int_{\Omega}e^{-2\lambda H_0(y)^2}\phi(y)^2\,dy,
\end{equation}
where $\lambda>0$, 
$g(x,t):=\lambda H_0(x)^2/(1-4 \lambda t)$ and $S_\lambda=1/4\lambda$. 
\end{proposition}
{\bf Proof.}
Similarly to \eqref{eq:3.10}, it follows that
\begin{equation*}
\begin{split}
 & \int_0^t\int_\Omega H(\nabla u)\nabla_\xi H(\nabla u)\nabla[e^{-2g}u]\,dy\,ds\\
 & \ge \int_0^t\int_\Omega e^{-2g}H(\nabla u)^2\,dy\,ds
 -2\int_0^t\int_\Omega e^{-2g}H(\nabla u)H(\nabla g)u\,dy\,ds\\
 & \ge (1-\mu)\int_0^t \int_\Omega e^{-2g}H(\nabla u)^2\,dy\,ds
 -\mu^{-1}\int_0^t\int_\Omega e^{-2g}u^2 H(\nabla g)^2 \,dy\,ds
\end{split}
\end{equation*}
for all $t\in(0,S_\lambda)$, where $\mu>0$. 
Then, similarly to \eqref{eq:3.11}, 
by \eqref{eq:3.1} we have
\begin{equation}
\label{eq:4.2}
\begin{split}
 & \frac{1}{2}\int_\Omega e^{-2g}u(y,s)^2\,dy\biggr|_{s=0}^{s=t}
 +(1-\mu)\int_0^t\int_\Omega e^{-2g}H(\nabla u)^2\,dy\,ds\\
 & \qquad\qquad
 \le \int_0^t\int_\Omega e^{-2g}u^2\left[\mu^{-1}H(\nabla g)^2-\partial_t g\right]\,dy\,ds
\end{split}
\end{equation}
for $t\in(0,S_\lambda)$. 
Setting $\mu=1$, 
by \eqref{eq:2.2} we see that
$$
\mu^{-1} H(\nabla g)^2-\partial_t g
=\frac{4\lambda^2}{(1-4\lambda t)^2}H_0(y)^2-\frac{4\lambda^2 H_0(y)^2}{(1-4\lambda t)^2}=0. 
$$
This together with \eqref{eq:4.2} implies that 
$$
 \frac{1}{2}\int_\Omega e^{-2g(y, t )}u(y, t )^2\,dy
 \le\frac{1}{2}\int_\Omega e^{-2g(y,0)}u(y,0)^2\,dy\\
$$
for all $t\in(0,S_\lambda)$. 
Thus \eqref{eq:4.1} holds. The proof is complete. 
$\Box$
\vspace{5pt}

Now we are ready to complete the proof of assertion~(ii) of Theorem~\ref{Theorem:1.2}.
\vspace{5pt}
\newline
{\bf Proof of assertion~(ii) of Theorem~\ref{Theorem:1.2}.} 
It suffices to consider the case $|\mu|({\bf R}^N)\not=0$. 
Let $\lambda>\Lambda>0$ and assume \eqref{eq:1.12}. 
Due to \eqref{eq:1.5}, we can assume, without loss of generality, that $\lambda>1>\Lambda>0$. 

By the Jordan decomposition theorem 
there exist two nonnegative Radon measure $\mu^+$ and $\mu^-$ such that $\mu=\mu^+-\mu^-$. 
Furthermore, we can find sequences $\{\mu_n^\pm\}\subset  C^\infty({\bf R}^N)$ 
such that 
\begin{equation}
\label{eq:4.3}
\lim_{n\to\infty}\int_{{\bf R}^N}\mu_n^\pm(x)\psi(x)\,dx=\int_{{\bf R}^N}\psi(x)\,d\mu^\pm (y)
\quad\mbox{for}\quad\psi\in C_0({\bf R}^N).
\end{equation}
For any $m=1,2,\dots$, let $\zeta_m\in C^\infty_0({\bf R}^N)$ be such that 
$$
\zeta_m=1\quad\mbox{in}\quad B_{H_0}(0,m/2),
\qquad
\mbox{supp}\,\zeta_m\subset B_{H_0}(0,m),
\qquad
0\leq \zeta_m \leq 1\quad\mbox{in}\quad{\bf R}^N. 
$$
Let $u_{m,n}$ be a solution of \eqref{eq:3.1} 
with $\Omega$ and $\phi$ replaced by $B_m:=B_{H_0}(0,m)$ and $\mu_{m,n}:=\zeta_m(\mu_n^+-\mu_n^-)  \in C^\infty_0(B_m)$, respectively. 
Setting $u_{m,n}(\cdot,t)=0$ outside $B_m$ for any $t>0$ and applying Propositions~\ref{Proposition:3.1} and \ref{Proposition:3.2}, 
we obtain 
\begin{equation}
\label{eq:4.4}
\begin{split}
 & 
\sup_{0<t< T_*}
\sup_{x\in{\bf R}^N}\int_{B_m \,\cap\,B_{H_0}(x,1)}e^{-h(y,t)}|u_{m,n}(y,t)|\,dy\\
 & \qquad
  +\sup_{0<t< T_*} t^{-\sigma'}\sup_{x\in{\bf R}^N}\int_0^t\int_{B_m\,\cap\,B_{H_0}(x,1)}e^{-h(y,s)}H(\nabla u_{m,n}(y,s))\,dy\,ds\\
 & \qquad\qquad
+\sup_{0<t< T_*}
t^{\frac{N}{4}}\sup_{x\in{\bf R}^N}\biggr(\int_{B_m\,\cap\,B_{H_0}(x,1)}|e^{-h(y,t)}u_{m,n}(y,t)|^2\,dy\biggr)^{\frac{1}{2}}\\
 & \le CI_{m,n}
 \quad\mbox{with}\quad
 I_{m,n}:=\sup_{x\in{\bf R}^N}\int_{B_m \,\cap\,B_{H_0}(x,1)}e^{-H_0(y)^2}|\mu_{m,n}(y)|\,dy
\end{split}
\end{equation}
for all $m$, $n=1,2,\dots$, where $T_*$ and $\sigma'$ are as in Propositions~\ref{Proposition:3.1} and \ref{Proposition:3.2}, respectively, 
and $h(x,t):= H_0(x)^2(1+t^\ell)$. 
Here and in what follows, $C$ denotes generic positive constant independent of $m$ and $n$.  

Let $1<\lambda'<\lambda$. Then we can find $t_*\in(0,T_*)$ such that $1+t_*^\ell<\lambda'$. 
Then, by \eqref{eq:4.4} we see that 
\begin{equation*}
\begin{split}
 & \int_{{\bf R}^N} e^{-2\lambda' H_0(y)^2}u_{m,n}(y,\tau)^2\,dy
 =\int_{B_m} e^{-2\lambda' H_0(y)^2}u_{m,n}(y,\tau)^2\,dy\\
 & \le C\sup_{x\in{\bf R}^N}\int_{B_m\,\cap\,B_{H_0}(x,1)} e^{-2h(y,\tau)}u_{m,n}(y,\tau)^2\,dy
\le C\tau^{-\frac{N}{2}}I_{m,n}^2
\end{split}
\end{equation*}
for $\tau\in(0,t_*]$. 
Applying Proposition~\ref{Proposition:4.1} with $u( y,  t)=u_{m,n}( y,  t+\tau)$  and $\lambda = \lambda'$, 
we obtain 
\begin{equation}
\label{eq:4.5}
\begin{split}
  \sup_{\tau< t < S_{\lambda'}+\tau}\int_{{\bf R}^N}e^{-\frac{2\lambda'H_0(y)^2}{1-4\lambda'( t -\tau)}}u_{m,n}(y, t )^2\,dy
 & \le C\int_{{\bf R}^N}e^{- 2\lambda' H_0(y)^2}u_{m,n}(y,\tau)^2\,dy\\
 & \le C\tau^{-\frac{N}{2}}I_{m,n}^2
\end{split}
\end{equation}
for all $\tau\in(0,t_*]$.
Taking a sufficiently small $t_*\in(0,T_*)$ if necessary, 
we see that 
\begin{equation*}
\begin{split}
 & 1+ t^\ell\le1+t_*^\ell<\lambda'<\lambda<\frac{\lambda}{1-4\lambda  t }
\quad\mbox{for}\quad 0< t <t_*,\\
 & \frac{1}{1-4\lambda'( t -t_*)}<\frac{1}{1-4\lambda  t }  - 2 \lambda t_*  \quad\mbox{for}\quad 0< t <S_\lambda.
\end{split}
\end{equation*}
Then, by \eqref{eq:4.4} and 
by applying the H\"older inequality to \eqref{eq:4.5} with $\tau=t_*$, 
we obtain 
\begin{equation}
\label{eq:4.6}
\begin{split}
\sup_{0< t <S_\lambda}\int_{{\bf R}^N}e^{-\lambda H_0(y)^2/(1-4\lambda  t )}|u_{m,n}(y, t )|\,dy
\le  C  I_{m,n}
\end{split}
\end{equation}
for $m$, $n=1,2,\dots$. 

On the other hand, 
it follows from \eqref{eq:1.12} with $\Lambda<1$ and \eqref{eq:4.3} that 
\begin{equation*}
\begin{split}
 & \limsup_{n\to\infty}I_{m,n}
\le \lim_{n\to\infty}\int_{B_m}e^{-H_0(y)^2}|\mu_{m,n}(y)|\,dy\\
 & = \lim_{n\to\infty}\int_{{\bf R}^N} e^{-H_0(y)^2}\zeta_m(y)  [\mu_n^+(y)+\mu_n^-(y)]\,dy
 =\int_{{\bf R}^N}e^{-H_0(y)^2} \zeta_m(y)  \, d|\mu|(y)\\
 & \leq \int_{{\bf R}^N} e^{-H_0(y)^2} \, d|\mu|(y)  \le C\sup_{x\in {\bf R}^N}\int_{B_{H_0}(x, 1/\sqrt{\Lambda})} e^{-\Lambda H_0(y)^2}\,d|\mu|(y)<\infty
\end{split}
\end{equation*}
for $m=1,2,\dots$. 
Since $|\mu|({\bf R}^N)\not=0$, by the diagonal argument 
we can find a sequence $\{n_m\}_{m=1}^\infty$ such that 
\begin{equation}
\label{eq:4.7}
\sup_m I_{m,n_m}
\le C\sup_{x \in {\bf R}^N}\int_{B_{H_0}(x,1/\sqrt{\Lambda})} e^{-\Lambda H_0(y)^2}d|\mu|(y)<\infty.
\end{equation}
This together with \eqref{eq:4.5} implies that 
$$
\sup_m \|u_{m,n_m}\|_{L^2(K)}<\infty
$$
for any compact set $K\subset{\bf R}^N\times(0,S_\lambda)$ and $m=1,2,\dots$. 
Then, by \eqref{eq:1.15} we apply the standard parabolic regularity theorems to the solution~$u_{m,n_m}$ 
and we can find $\alpha\in(0,1)$ such that  
\begin{equation}
\label{eq:4.8}
\sup_m\|u_{m,n_m}\|_{C^{1,\alpha;0,\alpha/2}(K)}<\infty 
\end{equation}
for any compact set $K\subset{\bf R}^N\times(0,S_\lambda)$. 
See \cite[Chapter~III]{LSU}. (See also \cite{DF1}, \cite{KM01}, \cite{KM02} and the last comment in \cite{CS}.)
Therefore, by the Ascoli-Arzel\`a theorem and the diagonal argument,  
taking a subsequence if necessary, 
we can find a function $u\in C^{1,\alpha;0,\alpha/2}({\bf R}^N\times(0,S_\lambda))$ such that 
\begin{equation}
 \lim_{m\to\infty}u_{m,n_m}(x,t)=u(x,t),
\quad
\lim_{m\to\infty}\nabla u_{m,n_m}(x,t)=\nabla u(x,t), 
\label{eq:4.9}
\end{equation}
uniformly on any compact set $K\subset{\bf R}^N\times(0,S_\lambda)$. 
Furthermore, due to \eqref{eq:4.4}, \eqref{eq:4.7} and \eqref{eq:4.9}, 
the Fatou lemma implies that 
$u\in L^1_{\rm loc}([0,S_\lambda):W^{1,1}(B_{H_0}(0,R)))$ for any $R>0$.
Moreover, by \eqref{eq:4.4} and \eqref{eq:4.7} we see that
\begin{equation}
\label{eq:4.10}
\begin{split}
 \int_0^\epsilon\int_{B_{H_0}(0,R)}
\left[|u_{m,n_m}|+H(\nabla u_{m,n_m})|\right]\,dy\,ds
\le C[\epsilon+\epsilon^{\sigma'}]
\end{split}
\end{equation}
for all sufficiently small $\epsilon>0$. 
Then, by \eqref{eq:4.8}, \eqref{eq:4.9} and \eqref{eq:4.10} 
we apply the Lebesgue dominated convergence theorem to obtain 
\begin{equation*}
\begin{split}
\lefteqn{\lim_{m\to\infty}\int_0^t\int_{{\bf R}^N}
 \left[-u_{m,n_m}\partial_t\varphi+H(\nabla u_{m,n_m})\nabla_\xi H(\nabla u_{m,n_m})\nabla\varphi\right]\,dy\,ds}\nonumber\\
 & = \lim_{m\to\infty}\int_0^\epsilon\int_{{\bf R}^N}
 \left[-u_{m,n_m}\partial_t\varphi+H(\nabla u_{m,n_m})\nabla_\xi H(\nabla u_{m,n_m})\nabla\varphi\right]\,dy\,ds\nonumber\\
 & \quad + \lim_{m\to\infty}\int_\epsilon^t\int_{{\bf R}^N}
 \left[-u_{m,n_m}\partial_t\varphi+H(\nabla u_{m,n_m})\nabla_\xi H(\nabla u_{m,n_m})\nabla\varphi\right]\,dy\,ds\nonumber\\
 & = \mathrm{O}(\epsilon+\epsilon^{\sigma'}) + \int_\epsilon^t\int_{{\bf R}^N}
 \left[-u\partial_t\varphi+H(\nabla u)\nabla_\xi H(\nabla u)\nabla\varphi\right]\,dy\,ds
\end{split}
\end{equation*}
for any $\varphi\in C_0^\infty({\bf R}^N\times[0,S_\lambda))$ and $0<t< S_\lambda$. 
Since $u \in L^1_{\rm loc}([0,S_\lambda) : W^{1,1}(B_{H_0}(0,R)))$ for any $R>0$ and $\epsilon$ is arbitrary, 
we deduce that
\begin{equation*}
 \begin{split}
 & \lim_{m\to\infty}\int_0^t\int_{{\bf R}^N}
 \left[-u_{m,n_m}\partial_t\varphi+H(\nabla u_{m,n_m})\nabla_\xi H(\nabla u_{m,n_m})\nabla\varphi\right]\,dy\,ds\\
 & = \int_0^t\int_{{\bf R}^N}
 \left[-u\partial_t\varphi+H(\nabla u)\nabla_\xi H(\nabla u)\nabla\varphi\right]\,dy\,ds.
 \end{split}
\end{equation*}
Furthermore, recalling the weak formulation of \eqref{eq:3.1} with $u = u_{m,n_m}$, 
we see that 
\begin{equation*}
\begin{split}
  & \lim_{m\to\infty}\int_{{\bf R}^N}u_{m,n_m}(y,t)\varphi(y,t)\,dy
+\int_0^t\int_{{\bf R}^N}\,\left[-u\partial_t\varphi+H(\nabla u)\nabla_\xi H(\nabla u)\nabla\varphi\right]\,dy\,ds\\
  & =\lim_{m\to\infty}\int_{{\bf R}^N}\varphi(y,0)\mu_{m,n_m}(y)\,dy
\end{split}
\end{equation*}
for all $\varphi\in C^\infty_0(B_{H_0}(0,R)\times[0,S_\lambda))$ and $0<t< S_\lambda$. 
This together with \eqref{eq:4.3} and \eqref{eq:4.9} implies that 
\begin{equation*}
\begin{split}
  & \int_{{\bf R}^N}u(y,t)\varphi(y,t)\,dy
+\int_0^t\int_{{\bf R}^N}\,\left[-u\partial_t\varphi+H(\nabla u)\nabla_\xi H(\nabla u)\nabla\varphi\right]\,dy\,ds\\
  & =\int_{{\bf R}^N}\varphi(y,0)\,d\mu(y) 
\end{split}
\end{equation*}
for all $0<t< S_\lambda$ and $\varphi\in C^\infty_0(B_{H_0}(0,R)\times[0, S_\lambda))$. 
In addition, by \eqref{eq:4.6}, \eqref{eq:4.7} and \eqref{eq:4.9} 
we apply the Fatou lemma again and see that $u$ satisfies \eqref{eq:1.13}. 
Therefore $u$ is the desired solution of \eqref{eq:1.6}. 
Thus assertion~(ii) of Theorem~\ref{Theorem:1.2} follows.  
$\Box$
\vspace{5pt}
\newline
{\bf Proof of assertion~(i) of Theorem~\ref{Theorem:1.2}.}
Let $u$ be a nonnegative solution of \eqref{eq:1.4} in ${\bf R}^N\times(0,T)$ for some $T>0$. 
By \eqref{eq:1.15}, applying the same argument as in \cite{IJK}, 
we can find a unique Radon measure $\mu$ in ${\bf R}^N$ such that 
$\mu$ satisfies \eqref{eq:1.11} and 
$$
\int_{{\bf R}^N}e^{-C|y|^2}d\mu(y)<\infty
$$
for some constant $C>0$. 
Since $H_0$ is an equivalent norm to the  Euclidean  norm $|\cdot|$ of ${\bf R}^N$, 
we obtain 
 $$
\int_{{\bf R}^N}e^{-C'H_0(y)^2}d\mu(y)<\infty. 
$$
Thus assertion~(i) of Theorem~\ref{Theorem:1.2} follows.
$\Box$
\vspace{5pt}

\noindent
{\bf Acknowledgements.} 
 G.A.~was partially supported by the Grant-in-Aid for Scientific Research (B)(No.~16H03946)
from Japan Society for the Promotion of Science. 
K.I.~would like to thank Professor Paolo Salani 
for his valuable suggestions. 
Furthermore, he also would like to thank Professor Juha Kinnunen 
for informing him of the papers \cite{KM01} and \cite{KM02}. 
He was partially supported 
by the Grant-in-Aid for Scientific Research (A)(No.~15H02058)
from Japan Society for the Promotion of Science. 
R.S.~was supported in part by Research Fellow 
of Japan Society for the Promotion of Science. 
\bibliographystyle{amsplain}

\begin{thebibliography}{10}

\bibitem{AD1} 
D. Andreucci and E. DiBenedetto, 
A new approach to initial traces in nonlinear filtration, 
Ann. Inst. H. Poincar\'e Anal. Non Lin\'eaire {\bf 7} (1990), 305--334. 

\bibitem{AD2}
D. Andreucci and E. DiBenedetto, 
On the Cauchy problem and initial traces for a class of evolution equations 
with strongly nonlinear sources, 
Ann. Scuola Norm. Sup. Pisa Cl. Sci. {\bf 18} (1991), 363--441.


\bibitem{Ar}
D. G. Aronson, 
Non-negative solutions of linear parabolic equations, 
Ann. Scuola Norm. Sup. Pisa {\bf 22} (1968), 607--694.

\bibitem{BCS}
D. Bao, S. S. Chern and Z. Shen, 
{\it An Introduction to Riemann--Finsler Geometry}, 
Springer, New York, 2000.


\bibitem{BaP} 
P. Baras and M. Pierre, 
Crit\`ere d'existence de solutions positives pour des \'equations semi-lin\'eaires non monotones, 
Ann. Inst. H. Poincar\'e Anal. Non Lin\'eaire {\bf 2} (1985), 185--212.

\bibitem{BP}
G. Bellettini and M. Paolini,  
Anisotropic motion by mean curvature in the context of Finsler geometry, 
Hokkaido Math. J. {\bf 25} (1996), 537--566. 

\bibitem{BCP}
P. B\'enilan, M.~G. Crandall and M. Pierre, 
Solutions of the porous medium equation in ${\bf R}^N$ 
under optimal conditions on initial values, 
Indiana Univ. Math. J. {\bf 33} (1984), 51--87.

\bibitem{BGS}
C. Bianchini, G. Ciraolo and P. Salani, 
An overdetermined problem for the anisotropic capacity, 
Calc. Var. Partial Differential Equations {\bf 55} (2016), Art. 84, 24 pp. 
  
\bibitem{Br}
H.~Br\'ezis, 
{\it Operateurs Maximaux Monotones et Semi-Groupes de Contractions
dans les Espaces de Hilbert}, Math Studies, Vol.5 North-Holland,
Amsterdam/New York, 1973. 
	
\bibitem{CS}
A. Cianchi and P. Salani, 
Overdetermined anisotropic elliptic problems, 
Math. Ann. {\bf 345} (2009), 859--881.

\bibitem{CFV} 
M. Cozzi, A. Farina and E. Valdinoci, 
Gradient bounds and rigidity results for singular, degenerate, anisotropic partial differential equations, 
Comm. Math. Phys. {\bf 331} (2014), 189--214. 

\bibitem{DPB}
F. Della Pietra and G. di Blasio, 
Blow-up solutions for some nonlinear elliptic equations involving a Finsler-Laplacian,  
Publ. Mat. {\bf 61} (2017), 213--238. 

\bibitem{DB02} 
E. DiBenedetto, 
{\it Degenerate parabolic equations}, 
Universitext, Springer-Verlag, New York, 1993.


\bibitem{DF1}
E. DiBenedetto and A. Friedman, 
H\"{o}lder estimates for nonlinear degenerate parabolic systems, 
J. Reine Angew. Math. {\bf 357} (1985), 1--22; 
J. Reine Angew. Math. {\bf 363} (1985), 217--220 (Addendum).


\bibitem{DH1} 
E. DiBenedetto and M. A. Herrero,  
On the Cauchy problem and initial traces for a degenerate parabolic equation, 
Trans. Amer. Math. Soc. {\bf 314} (1989), 187--224. 

\bibitem{DH2} 
E. DiBenedetto and M. A. Herrero, 
Nonnegative solutions of the evolution $p$-Laplacian equation. Initial traces 
and Cauchy problem when $1<p<2$, 
Arch. Rational Mech. Anal. {\bf 111} (1990), 225--290.

\bibitem{EG} 
L. C. Evans and R. F.  Gariepy, 
{\it Measure Theory and Fine Properties of Functions}, 
Studies in Advanced Mathematics, CRC Press, Boca Raton, FL, 1992.

\bibitem{FK}
V. Ferone and B. Kawohl, 
Remarks on a Finsler-Laplacian, 
Proc. Amer. Math. Soc. {\bf 137} (2009), 247--253. 

\bibitem{GU}
Y. Giga and N. Umeda, 
On instant blow-up for semilinear heat equations with growing initial data, 
Methods Appl. Anal. {\bf 15} (2008), 185--195.

\bibitem{HP}
M.~A. Herrero and M. Pierre, 
The Cauchy problem for $u_t=\Delta u^m$ when $0<m<1$, 
Trans. Amer. Math. Soc. {\bf 291} (1985), 145--158.


\bibitem{HI02} 
K. Hisa and K. Ishige, 
Solvability of the heat equation with a nonlinear boundary condition, 
preprint (arXiv:1704.07992). 

\bibitem{I0}
K. Ishige, 
On the existence of solutions of the Cauchy problem 
for porous medium equations with Radon measure as initial data,
Discrete Contin. Dynam. System {\bf 1} (1995), 521--546.
	
\bibitem{I}
K. Ishige, 
On the existence of solutions of the Cauchy problem 
for a doubly nonlinear parabolic equation, 
SIAM J. Math. Anal. {\bf 27} (1996), 1235--1260. 

\bibitem{I1}
K. Ishige, 
On the existence of solutions of the Cauchy problem 
for a quasi-linear parabolic equation with unbounded initial data, 
Adv. Math. Sci. Appl. {\bf 9} (1999), 263--289.

\bibitem{IJK}
K. Ishige and J. Kinnunen, 
Initial trace for a doubly nonlinear parabolic equation, 
J. Evol. Equ. {\bf 11} (2011), 943--957. 

\bibitem{IS1}
K. Ishige and R. Sato, 
Heat equation with a nonlinear boundary condition
and uniformly local $L^r$ spaces, 
Discrete Contin. Dyn. Syst. {\bf 36} (2016), 2627--2652.

\bibitem{IS2}
K. Ishige and R. Sato, 
Heat equation with a nonlinear boundary condition and growing initial data, 
Differential Integral Equations {\bf 30} (2017), 481--504.

\bibitem{KM01}
T. Kuusi and G. Mingione, 
Gradient regularity for nonlinear parabolic equations, 
Ann. Sc. Norm. Super. Pisa Cl. Sci. {\bf 12} (2013), 755--822. 

\bibitem{KM02}
T. Kuusi and G. Mingione, 
The Wolff gradient bound for degenerate parabolic equations, 
J. Eur. Math. Soc. {\bf 16} (2014), 835--892.

\bibitem{LSU}
O. A. Lady\v{z}enskaja, V. A. Solonnikov and N. N. Ural'ceva,
{\it Linear and Quasi-linear Equations of Parabolic Type},
American Mathematical Society Translations, vol. 23, American Mathematical Society,
Providence, RI, 1968.

\bibitem{OS1}
S.-I. Ohta and K.-T. Sturm, 
Heat flow on Finsler manifolds, 
Comm. Pure Appl. Math. {\bf 62} (2009), 1386--1433. 

\bibitem{OS2}
S.-I. Ohta and K.-T. Sturm, 
Bochner-Weitzenb\"ock formula and Li-Yau estimates on Finsler manifolds, 
Adv. Math. {\bf 252} (2014), 429--448. 

\bibitem{S}
R. Schneider, 
{\it Convex Bodies: The Brunn-Minkowski Theory}, 
Cambridge University Press, Cambridge (1993). 

\bibitem{T}
A. Tachikawa, 
A partial regularity result for harmonic maps into a Finsler manifold, 
Calc. Var. Partial Differential Equations {\bf 16} (2003), 225--226. 

\bibitem{W} 
D.~V. Widder, 
Positive temperatures on an infinite rod, 
Trans. Amer. Math. Soc. {\bf 55} (1944), 85--95. 

\bibitem{X}
Q. Xia, 
On the first eigencone for the Finsler Laplacian, 
Bull. Aust. Math. Soc. {\bf 94} (2016), 316--325. 

\bibitem{Y} 
K. Yamada, 
On viscous conservation laws with growing initial data, 
Differential Integral Equations {\bf 18} (2005), 841--854. 

\end{thebibliography}

\end{document}